\documentclass[11pt,leqno]{article} 

\usepackage{fullpage}

\usepackage{amssymb} 
\usepackage{amsmath} 
\date{}
\title{Classification of linearly compact simple Nambu-Poisson  
 algebras} 
\author{{\sc Nicoletta Cantarini}\thanks{Dipartimento di Matematica, Universit\`a di Bologna, Bologna, Italy}
\and
{\sc Victor G.\ Kac}\thanks{Department of Mathematics, MIT, Cambridge,
Massachusetts 02139, USA}}

\newtheorem{theorem}{Theorem}[section] 
\newtheorem{lemma}[theorem]{Lemma} 
\newtheorem{corollary}[theorem]{Corollary} 
\newtheorem{proposition}[theorem]{Proposition} 
\newtheorem{definition}[theorem]{Definition} 
\newtheorem{remark}[theorem]{Remark}

\newtheorem{example}[theorem]{Example}

\def\Z{\mathbb{Z}}

\def\F{\mathbb{F}}

\def\0{\bar{0}}
\def\1{\bar{1}}

\numberwithin{equation}{section}

\makeatletter

\def\enumerate{%
  \ifnum \@enumdepth >\thr@@\@toodeep\else
    \advance\@enumdepth\@ne
    \edef\@enumctr{enum\romannumeral\the\@enumdepth}%
      \list
        {\csname label\@enumctr\endcsname}%
        {\usecounter\@enumctr
          \addtolength{\leftmargin}{-\leftmargin}
          \settowidth{\labelwidth}{(99)}
          \itemindent = \labelwidth
          \addtolength{\itemindent}{\labelsep}
        \listparindent=1em      
          \def\makelabel##1{{##1}\hfill}
          }%
  \fi}

\makeatother

\begin{document} 
\maketitle 
\begin{abstract} We introduce  the notion of universal odd generalized Poisson superalgebra associated to an associative algebra $A$,
 by generalizing a construction made in \cite{DeSK}. By making use of this notion we give a complete classification of simple linearly compact (generalized) $n$-Nambu-Poisson algebras over an algebraically closed field of characteristic zero.
\end{abstract} 
\section*{Introduction}
In 1973 Y.\ Nambu proposed a generalization of Hamiltonian mechanics, based on the notion of $n$-ary bracket in place
of the usual binary Poisson bracket \cite{N}. Nambu dynamics is described by the flow, given by a system of ordinary differential equations which involves $n-1$ Hamiltonians:
\begin{equation}
\frac{du}{dt}=\{u,h_1,\dots, h_{n-1}\}.
\label{0.1}
\end{equation}
The (only) example, proposed by Nambu is the following $n$-ary bracket on the space of functions in $N\geq n$ variables:
\begin{equation}
\{f_1, \dots, f_n\}=\det\left(\frac{\partial f_i}{\partial x_j}\right)^n_{i,j=1}.
\label{0.2}
\end{equation}
He pointed out that this $n$-ary bracket satisfies the following axioms, similar to that of a Poisson bracket:
$$(\mbox{Leibniz rule})~~~~~\{f_1,\dots, f_i\tilde{f}_i,\dots, f_n\}=
f_i\{f_1,\dots, \tilde{f}_i,\dots, f_n\}+\tilde{f}_i\{f_1,\dots, f_i,\dots, f_n\};$$
$$(\mbox{skewsymmetry})~~~~~\{f_{\sigma(1)},\dots, f_{\sigma(n)}\}=(sign \sigma)
\{f_1,\dots, f_n\}.$$
Twelve years later this example was rediscovered by F.\ T.\ Filippov in his theory of $n$-Lie algebras which
is a natural generalization of ordinary (binary) Lie algebras \cite{F}. Namely, an $n$-Lie algebra is a vector space with $n$-ary bracket $[a_1, \dots, a_n]$, which is skewsymmetric (as above) and satisfies the following Filippov-Jacobi identity:
\begin{equation}
[a_1,\dots , a_{n-1}, [b_1,\dots , b_n]] = [[a_1, \dots , a_{n-1}, b_1], b_2, \dots , b_n] + [b_1, [a_1,\dots , a_{n-1}, b_2], b_3, \dots, b_n]
+\dots
\label{0.3}
\end{equation}
$$ + [b_1, \dots , b_{n-1}, [a_1,\dots , a_{n-1}, b_n]].$$
In particular, Filippov proved that the Nambu bracket (\ref{0.2}) satisfies the Filippov-Jacobi identity.

Following Takhtajan \cite{T}, we call an {\em $n$-Nambu-Poisson algebra} a unital commutative associative algebra ${\cal N}$, endowed with an $n$-ary bracket, satisfying the Leibniz rule, skew-symmetry and Filippov-Jacobi identity. Of course for $n=2$ this is the definition of a Poisson algebra.

In \cite{CantaK4} we classified simple linearly compact $n$-Lie algebras with $n>2$ over a field $\F$ of characteristic $0$. The classification is based on a bijective correspondence between $n$-Lie algebras and pairs $(L, \mu)$, where $L$ is a $\Z$-graded Lie superalgebra of the form $L=\oplus_{j=-1}^{n-1}L_j$ satisfying certain additional properties, and $L_{n-1}=\F\mu$, thereby reducing it to the known classification of simple linearly compact Lie superalgebras and their $\Z$-gradings  \cite{K3}, \cite{CantaK}. For this construction we used the universal $\Z$-graded Lie superalgebra, associated to a vector superspace.

In the present paper we use an analogous correspondence between linearly compact $n$-Nambu-Poisson algebras and certain "good" pairs 
$({\cal P}, \mu)$, where ${\cal P}$ is a $\Z_+$-graded odd Poisson superalgebra ${\cal P}=\oplus_{j\geq -1}{\cal P}_j$ and 
$\mu\in {\cal P}_{n-1}$  is an element of parity $n$ mod $2$. For this construction we use the universal $\Z$-graded odd Poisson
superalgebra, associated to an associative algebra, considered in \cite{DeSK}.
As a result, using the classification of simple linearly compact odd Poisson superalgebras \cite{CantaK3}, we obtain the following theorem.
\begin{theorem}\label{T01}
For $n>2$, any simple linearly compact $n$-Nambu-Poisson algebra is isomorphic to the algebra $\F[[x_1, \dots, x_n]]$ with the $n$-ary bracket (\ref{0.2}). 
\end{theorem}
Note the sharp difference with the Poisson case, when each algebra $\F[[p_1, \dots, p_n, q_1, \dots, q_n]]$ carries a 
Poisson bracket
\begin{equation}
\{f,g\}_P=\sum_{i=1}^n(\frac{\partial f}{\partial p_i}\frac{\partial g}{\partial q_i}-
\frac{\partial f}{\partial q_i}\frac{\partial g}{\partial p_i}),
\label{0.4}
\end{equation}
making it a simple linearly compact Poisson algebra (and these are all, up to isomorphism \cite{CantaK2}).

In the present paper we treat also the case of a generalized $n$-Nambu-Poisson bracket, which is an $n$-ary analogue of the 
generalized Poisson bracket, called also the Lagrange's bracket. For the latter bracket the Leibniz rule is modified by adding 
an extra term:
$$\{a, bc\}=\{a,b\}c+\{a,c\}b-\{a,1\}bc.$$
In order to treat this case along similar lines, we construct the universal $\Z$-graded generalized odd Poisson superalgebra,
associated to an associative algebra, which is a generalization of the construction in \cite{DeSK}. Our main result in this direction is the following theorem, which uses the classification of simple linearly compact odd generalized Poisson superalgebras \cite{CantaK3}.
\begin{theorem}\label{T02}
For $n>2$, any simple linearly compact generalized $n$-Nambu-Poisson algebra is gauge equivalent (see Remark \ref{gaugeNambu} for the definition) either to the Nambu $n$-algebra from 
Theorem \ref{T01} or to the Dzhumadildaev $n$-algebra \cite{D}, which is $\F[[x_1, \dots, x_{n-1}]]$ with the 
$n$-ary bracket 
\begin{equation}
\{f_1, \dots, f_n\}=\det\left(
\begin{array}{ccc}
f_1 & \dots & f_n\\
\frac{\partial f_1}{\partial x_1} & \dots & \frac{\partial f_n}{\partial x_1}\\
\dots & \dots & \dots\\
\frac{\partial f_1}{\partial x_{n-1}} & \dots & \frac{\partial f_n}{\partial x_{n-1}}
\end{array}\right).
\label{0.5}
\end{equation}
\end{theorem}
Note again the sharp difference with the generalized Poisson case, when each algebra $\F[[p_1, \dots, p_n,$ $q_1, \dots, q_n,t]]$
carries a Lagrange bracket
\begin{equation}
\{f,g\}_L=\{f,g\}_P+(2-E)f\frac{\partial g}{\partial t}-\frac{\partial f}{\partial t}(2-E)g,
\label{0.6}
\end{equation}
where $\{f,g\}_P$ is given by (\ref{0.4}) and $E=\sum_{i=1}^n(p_i\frac{\partial}{\partial p_i}+q_i\frac{\partial}{\partial q_i})$,
making it a simple linearly compact generalized Poisson algebra (and those, along with (\ref{0.4}), are all, up to gauge equivalence).

Throughout the paper our base field $\F$ has characteristic $0$ and is algebraically closed.

\section{Nambu-Poisson algebras}
\begin{definition}\label{NambuPoisson} A generalized $n$-Nambu-Poisson algebra  is a triple $({\cal N},\{\cdot,\dots,\cdot\},\cdot)$
such that 
\begin{itemize}
\item[-] $({\cal N}, \cdot)$ is a unital
associative commutative algebra;
\item[-] 
$({\cal N},\{\cdot,\dots,\cdot\})$ is an $n$-Lie algebra; 
\item[-] the following generalized Leibniz rule
holds:
\begin{equation}
\{a_1,\dots, a_{n-1},bc\}=\{a_1, \dots, a_{n-1},b\}c+b\{a_1,\dots,a_{n-1},c\}-\{a_1, \dots, a_{n-1},1\}bc.
\label{nambuLeibniz}
\end{equation}
\end{itemize}
If $\{a_1, \dots, a_{n-1},1\}=0$, then (\ref{nambuLeibniz}) is the usual Leibniz rule and $({\cal N},\{\cdot,\dots,\cdot\},\cdot)$
is called simply $n$-Nambu-Poisson algebra.
\end{definition}

For $n=2$ Definition \ref{NambuPoisson}
is the definition of a generalized Poisson algebra. Simple linearly compact generalized Poisson (super)algebras were classified in
\cite[Corollary 7.1]{CantaK2}.

\begin{example}\label{Nambu}\em Let ${\cal N}=\F[[x_1,\dots, x_n]]$ with the usual commutative associative product and $n$-ary bracket
defined, for $f_1,\dots, f_n\in {\cal N}$, by: 
$$\{f_1, \dots, f_n\}=\det\left(\begin{array}{ccc}
D_1(f_1) & \dots & D_1(f_n)\\
\dots & \dots & \dots\\
D_n(f_1) & \dots & D_n(f_n)
\end{array}
\right)$$
where $D_i=\frac{\partial}{\partial x_i}$, $i=1,\dots, n$.
Then ${\cal N}$ is an $n$-Nambu-Poisson algebra, introduced by Nambu \cite{N}, that we will call the $n$-Nambu algebra (cf.\ \cite{N}, \cite{F},  \cite{CantaK4}). 
\end{example}

\begin{example}\label{Dzhuma}\em Let ${\cal N}=\F[[x_1,\dots, x_{n-1}]]$ with the usual commutative associative product and $n$-ary bracket
defined, for $f_1,\dots, f_n\in {\cal N}$, by 
$$\{f_1, \dots, f_n\}=\det\left(\begin{array}{ccc}
f_1 & \dots & f_n\\
D_1(f_1) & \dots & D_1(f_n)\\
\dots & \dots & \dots\\
D_{n-1}(f_1) & \dots & D_{n-1}(f_n)
\end{array}
\right)$$
where $D_i=\frac{\partial}{\partial x_i}$, $i=1,\dots, n-1$.
Then ${\cal N}$ is a generalized Nambu-Poisson algebra that we will call the
$n$-Dzhumadildaev algebra (cf.\ \cite{D}, \cite{CantaK4}). 
\end{example}


\begin{remark}\label{gaugeNambu}\em  Let $N=({\cal N},\{\cdot,\dots,\cdot\},\cdot)$ be a generalized $n$-Nambu-Poisson algebra. For any invertible
element $\varphi\in{\cal N}$ define the following bracket on ${\cal N}$:
\begin{equation}
\{f_1,\dots,f_n\}^{\varphi}=\varphi^{-1}\{\varphi f_1,\dots, \varphi f_n\}.
\label{gaugequivnambu}
\end{equation}
Then $N^{\varphi}=({\cal N},\{\cdot,\dots,\cdot\}^{\varphi},\cdot)$ is another generalized $n$-Nambu-Poisson algebra. Indeed,
the skew-symmetry of the bracket is straightforward and the Filippov-Jacobi identity for the bracket $\{\cdot, \dots, \cdot\}^\varphi$
easily follows from the Filippov-Jacobi identity for the bracket $\{\cdot, \dots, \cdot\}$.
Let us check that $\{\cdot, \dots, \cdot\}^\varphi$ satisfies the generalized Leibniz rule. We have:

\bigskip

\noindent
$\{f_1, \dots, f_{n-1}, gh\}^\varphi=\varphi^{-1}\{\varphi f_1, \dots, \varphi f_{n-1}, \varphi gh\}
=\varphi^{-1}(\{\varphi f_1, \dots, \varphi f_{n-1}, \varphi g\}h$

\bigskip

\noindent
$+\varphi g\{\varphi f_1, \dots, \varphi f_{n-1}, h\}-\{\varphi f_1, \dots, \varphi f_{n-1},1\} \varphi gh)
=\{f_1, \dots, f_{n-1}, g\}^\varphi h$

\bigskip

\noindent
$+g\{\varphi f_1, \dots, \varphi f_{n-1}, h\}-
\{\varphi f_1, \dots, \varphi f_{n-1},1\} gh=\{f_1, \dots, f_{n-1}, g\}^\varphi h$

\bigskip

\noindent
$+g\{\varphi f_1, \dots, \varphi f_{n-1}, h\}-
\{\varphi f_1, \dots, \varphi f_{n-1},1\} gh+g\{f_1, \dots, f_{n-1}, h\}^\varphi$

\bigskip

\noindent
$-
\varphi^{-1}g\{\varphi f_1, \dots, \varphi f_{n-1}, \varphi h\}=\{f_1, \dots, f_{n-1}, g\}^\varphi h+g\{\varphi f_1, \dots, \varphi f_{n-1}, h\}$

\bigskip

\noindent
$-\{\varphi f_1, \dots, \varphi f_{n-1},1\} gh+g\{f_1, \dots, f_{n-1}, h\}^\varphi-
g\{\varphi f_1, \dots, \varphi f_{n-1}, h\}$

\bigskip

\noindent
$-\varphi^{-1}gh\{\varphi f_1, \dots, \varphi f_{n-1}, \varphi\}+
\{\varphi f_1, \dots, \varphi f_{n-1},1\} gh$

\bigskip

\noindent
$
=\{f_1, \dots, f_{n-1}, g\}^\varphi h+g\{f_1, \dots, f_{n-1}, h\}^\varphi-\{f_1, \dots, f_{n-1}, 1\}^{\varphi}gh.$

\bigskip

\noindent
We shall say that the generalized Nambu-Poisson algebras $N$ and $N^{\varphi}$  are {\it gauge equivalent}.

\medskip
\end{remark}

\section{Odd generalized Poisson superalgebras}
\begin{definition}\label{OGP}
An odd generalized Poisson superalgebra $({\cal P}, [\cdot,\cdot], \wedge)$ is a triple
such that 
\begin{itemize}
\item[-] $({\cal P}, \wedge)$ is a unital
associative commutative superalgebra with parity $p$;
\item[-] 
$(\Pi{\cal P}, [\cdot,\cdot])$ is a Lie superalgebra (here $\Pi{\cal P}$ denotes the space ${\cal P}$ with parity $\bar{p}=p+\bar{1}$); 
\item[-] the following generalized odd Leibniz rule
holds:
\begin{equation}
[a,b\wedge c]=[a,b]\wedge c+(-1)^{(p(a)+1)p(b)}b\wedge [a,c]+(-1)^{p(a)+1}D(a)\wedge b\wedge c,
\label{oddleibniz}
\end{equation}
\end{itemize}
where $D(a)=[1,a]$. If $D=0$, then relation (\ref{oddleibniz}) becomes the odd Leibniz rule; in this case
$({\cal P}, [\cdot,\cdot], \wedge)$ is called an odd Poisson superalgebra (or Gerstenhaber superalgebra). 
Note that $D$ is an odd derivation of the associative product and of the Lie superalgebra bracket.
\end{definition}

\begin{example}\label{PO(n,n)}\em
Consider the
commutative associative superalgebra 
${\cal O}(m,n)=\Lambda(n)[[x_1,\dots, x_m]]$, 
where $\Lambda(n)$
denotes the Grassmann algebra over $\F$ on $n$ anti-commuting indeterminates
$\xi_1, \dots, \xi_n$, and the superalgebra parity is defined by
$p(x_i)=\bar{0}$, $p(\xi_j)=\bar{1}$.

Set $m=n$ and define the following bracket, known as the Buttin
bracket, on ${\cal O}(n,n)$
($f,g\in {\cal O}(n,n)$):
\begin{equation}[f,g]_{HO}=\sum_{i=1}^n(\frac{\partial f}{\partial x_i}\frac{\partial
g}{\partial \xi_i}+(-1)^{p(f)}\frac{\partial f}{\partial \xi_i}\frac{\partial
g}{\partial x_i}).
\label{buttinbracket}
\end{equation}
Then ${\cal O}(n,n)$ with this bracket is an odd Poisson superalgebra,
 which we denote  by $PO(n,n)$.
\end{example}

\begin{example}\label{PO(n,n+1)}\em
Consider the associative superalgebra ${\cal O}(n,n+1)$
with even indeterminates
$x_1, \dots, x_n$
 and odd indeterminates $\xi_1,\dots, \xi_n, \xi_{n+1}=\tau$.
Define on ${\cal O}(n, n+1)$ the following bracket
($f,g\in {\cal O}(n,n+1)$):
\begin{equation}
[f,g]_{KO}=[f,g]_{HO}+(E-2)(f)\frac{\partial g}{\partial \tau}
+(-1)^{p(f)}\frac{\partial
f}{\partial \tau}(E-2)(g),
\label{Leibniz!!}
\end{equation}
where $[\cdot,\cdot]_{HO}$ is the Buttin bracket (\ref{buttinbracket})
and $E=\sum_{i=1}^n(x_i\frac{\partial}{\partial
x_i}+\xi_i
\frac{\partial}{\partial \xi_i})$ is the Euler operator.
Then ${\cal O}(n,n+1)$ with bracket $[\cdot,\cdot]_{KO}$
is an odd  generalized Poisson superalgebra
with $D=-2\frac{\partial}{\partial \tau}$ \cite[Remark 4.1]{CantaK},
which we denote  by $PO(n,n+1)$.   
\end{example}

\begin{remark}\em Let $P=({\cal P},[\cdot,\cdot],\cdot)$ be an odd generalized Poisson superalgebra. For any invertible
element $\varphi\in{\cal P}$, such that $p(\varphi)=\bar{0}$ and $[\varphi,\varphi]=0$, define the following bracket on $P$:
\begin{equation}
[a,b]^{\varphi}=\varphi^{-1}[\varphi a, \varphi b].
\label{gaugequiv}
\end{equation}
Then $P^{\varphi}=({\cal P},[\cdot,\cdot]^{\varphi},\cdot)$ is another odd generalized Poisson superalgebra, with derivation 
$$D_{\varphi}(a)=
[1, a]^{\varphi}=[\varphi, a]-D(\varphi)a.$$ The odd generalized Poisson superalgebras $P$ and $P^{\varphi}$ are called
{\it gauge equivalent} (cf.\ \cite[Example 3.4]{CantaK3}). Note that the associative products in $P$ and $P^{\varphi}$ are the same. 
\end{remark}

\begin{theorem}\cite[Corollary 9.2]{CantaK3}\label{rigid}
\begin{itemize}
\item[a)] Any simple linearly compact odd generalized Poisson superalgebra is gauge
equivalent to $PO(n,n)$ or $PO(n,n+ 1)$.
\item[b)] Any simple linearly compact odd Poisson superalgebra is isomorphic to
$PO(n,n)$.
\end{itemize}
\end{theorem}

\begin{definition}\label{graded}
A $\Z$-graded (resp.\ $\Z_+$-graded)
 odd generalized  Poisson superalgebra is an odd generalized Poisson superalgebra
$({\cal P}, [\cdot,\cdot], \wedge)$ such that $(\Pi{\cal P},[\cdot,\cdot])$ is a $\Z$-graded Lie superalgebra:
 $\Pi{\cal P}=\oplus_{j\in\Z}{\cal P}_j$ (resp.\ a $\Z$-graded Lie superalgebra of depth 1: $\Pi{\cal P}=\oplus_{j\geq -1}{\cal P}_j$) 
 and $({\cal P}, \wedge)$ is a $\Z$-graded commutative associative superalgebra:
${\cal P}=\oplus_{k\in\Z}{\cal Q}_{k}$ (resp.\ a $\Z_+$-graded commutative associative superalgebra:
${\cal P}=\oplus_{k\in \Z_+}{\cal Q}_{k}$)  such that ${\cal P}_{j}=\Pi{\cal Q}_{j+1}$.
\end{definition}

\begin{example}\label{GradedGeneralized}\em
Let us consider the odd Poisson superalgebra $PO(n,n)$ (resp.\ $PO(n,n+1)$). Set
$\deg x_i=0$ and $\deg \xi_i=1$ for every $i=1, \dots, n$ (resp.\ $\deg x_i=0$, $\deg \xi_i=1$ for every $i=1, \dots, n$ and $\deg \tau=1$). Then $PO(n,n)$ (resp.\ $PO(n,n+1)$) becomes a $\Z_+$ graded odd (resp.\ generalized) Poisson
superalgebra with 
$${\cal Q}_j=\{f\in {\cal O}(n,n) ~|~\deg(f)=j\}$$
and
$${\cal P}_j=\{f\in {\cal O}(n,n)~|~ \deg(f)=j+1\}.$$
We will call this grading a grading of type $(0,\dots,0|1,\dots,1)$ (resp.\ $(0,\dots,0|1,\dots,1,1)$).
We thus have, for ${\cal P}=PO(n,n)$:
$$\Pi{\cal P}_{-1}={\cal Q}_0=\F[[x_1, \dots, x_n]]$$
and, for $j\geq 0$,
$$\Pi{\cal P}_j={\cal Q}_{j+1}=\langle \xi_{i_1}\dots\xi_{i_{j+1}} ~|~ 1\leq i_1< \dots <i_{j+1}\leq n\rangle\otimes\F[[x_1, \dots, x_n]].$$
Similarly, for ${\cal P}=PO(n,n+1)$, we have:
$${\cal P}_{-1}={\cal Q}_0=\F[[x_1, \dots, x_n]]$$
$${\cal P}_j={\cal Q}_{j+1}=\langle \xi_{i_1}\dots\xi_{i_{j+1}} ~|~ 1\leq i_1< \dots <i_{j+1}\leq n+1\rangle\otimes\F[[x_1, \dots, x_n]].$$
\end{example}

\begin{remark}\label{gradings}\em
From the properties of the $\Z$-gradings of the Lie superalgebras $HO(n,n)$ and $KO(n,$ $n+1)$ (see, for example, \cite{K3}),
one can deduce that the grading of type $(0,\dots,0|1,\dots,1)$ (resp.\ $(0,\dots,0|$ $1,\dots,1,1)$) is, up to isomorphisms,
the only $\Z_+$-grading of ${\cal P}=PO(n,n)$ (resp.\ ${\cal P}=PO(n,n+1)$) such that ${\cal P}_{-1}$ is completely odd.
\end{remark}

\begin{remark}\label{gaugegradings}\em
Let ${\cal P}=PO(n,n)$ or ${\cal P}=PO(n,n+1)$ and let ${\cal P}^{\varphi}$ be an odd generalized Poisson superalgebra gauge equivalent to ${\cal P}$.
Then
the grading of type $(0,\dots,0|1,\dots,1)$ (resp.\ $(0,\dots,0|$ $1,\dots,1,1)$) is, up to isomorphisms,
the only $\Z_+$-grading of ${\cal P}^{\varphi}$ such that ${\cal P}^{\varphi}_{-1}$ is completely odd.
Indeed, 
let ${\cal P}^{\varphi}=\oplus_{k\in \Z_+}{\cal Q}^{\varphi}_{k}=\oplus_{j\geq -1}{\cal P}^{\varphi}_{j}$, with  ${\cal P}^{\varphi}_{j}=\Pi{\cal Q^{\varphi}}_{j+1}$
a $\Z_+$-grading of ${\cal P}^{\varphi}$. Suppose that, $x_i\in {\cal Q}^{\varphi}_{k}$ and $\xi_i\in {\cal Q}^{\varphi}_{j}$ for some $1\leq i\leq n$ and some
$k,j\in\Z_+$. Then 
\begin{equation}
[x_i,\xi_i]^{\varphi}\in \Pi {\cal P}^{\varphi}_{k+j-2}={\cal Q}^{\varphi}_{k+j-1}.
\label{comm.rel}
\end{equation}
 On the other hand, by (\ref{gaugequiv}),
we have:
$$[x_i,\xi_i]^{\varphi}=\varphi^{-1}[\varphi x_i, \varphi\xi_i]=
\varphi^{-1}([\varphi x_i,\varphi]\xi_i+\varphi[\varphi x_i,\xi_i]-D(\varphi x_i)\varphi\xi_i)=$$
$$=[x_i,\varphi]\xi_i-D(\varphi) x_i\xi_i+[\varphi x_i, \xi_i]-
D(\varphi x_i)\xi_i=\frac{\partial\varphi}{\partial \xi_i}\xi_i+\frac{1}{2}x_i\xi_iD(\varphi)-D(\varphi)x_i\xi_i+
\frac{\partial \varphi}{\partial x_i}x_i+$$
$$+\varphi+\frac{1}{2}D(\varphi)x_i\xi_i-D(\varphi)x_i\xi_i=\frac{\partial \varphi}{\partial \xi_i}\xi_i-D(\varphi)x_i\xi_i+\frac{\partial \varphi}{\partial x_i}x_i+\varphi,$$
where $D=0$ if ${\cal P}=PO(n,n)$ and $D=-2\frac{\partial}{\partial \tau}$ if ${\cal P}=PO(n,n+1)$.
Note that $[x_i,\xi_i]^{\varphi}$ is invertible since $\varphi$ is invertible and, by (\ref{comm.rel}), it is homogeneous, hence
$k+j=1$, i.e., either $k=0$ and $j=1$ or $k=1$ and $j=0$. It follows that the only $\Z_+$grading of ${\cal P}^{\varphi}$ 
 such that ${\cal P}^{\varphi}_{-1}$ is completely odd is the grading of type $(0,\dots,0|1,\dots,1)$. We can thus simply denote the graded components 
of ${\cal P}^{\varphi}$ with respect to this grading by ${\cal P}_{j}=\Pi{\cal Q}_{j+1}$. 

Now let $a\in {\cal Q}_i=\Pi {\cal P}_{i-1}$ and $b\in {\cal Q}_k=\Pi {\cal P}_{k-1}$. 
We have: $[a,b]^{\varphi}=[a,\varphi]b+[\varphi a,b]+(-1)^{p(a)+1}(D(\varphi)ab+D(\varphi a)b)$. 
Suppose that $\varphi=\sum_{j\geq 0}\varphi_j$ with $\varphi_j\in {\cal Q}_j$. Then one can show, using the fact that
 $[a,b]^{\varphi}\in \Pi {\cal P}_{i+k-2}={\cal Q}_{i+k-1}$, that $[a,b]^{\varphi}=[a,b]^{\varphi_0}$.
It follows that when dealing with the $\Z_+$-graded odd generalized Poisson superalgebras ${\cal P}^{\varphi}$ we can always assume
$\varphi\in{\cal Q}_0$.
\end{remark}

\section{The universal odd generalized Poisson superalgebra}

\begin{definition}
Let $A$ be a unital commutative associative superalgebra with parity $p$. A linear map $X: A \rightarrow A$ is called a generalized derivation of $A$ if it satisfies the generalized Leibniz rule:
\begin{equation}\label{genleibniz}
X(bc) = X(b)c + (-1)^{p(b)p(c)}X(c)b - X(1)bc.
\end{equation}
We denote by $GDer(A)$ the set of generalized derivations of $A$. If $X(1) = 0$, relation (\ref{genleibniz}) becomes the usual Leibniz rule and $X$ is called a derivation. We denote by $Der(A)$ the set of derivations of $A$. 
\end{definition}

\begin{proposition}
The set $GDer(A)$ is a subalgebra of the Lie superalgebra $End(A)$.
\end{proposition}
{\bf Proof.} This follows by direct computations.
\hfill$\Box$

\medskip

Our construction of the universal odd generalized Poisson superalgebra is inspired by the one of the universal odd Poisson superalgebra explained in \cite{DeSK}. The universal odd Poisson superalgebra associated to $A$ is the full prolongation of the subalgebra $Der(A)$ of the Lie superalgebra $End(A)$ (the definitions will be given below). In this section we generalize this construction when $Der(A)$ is replaced by the subalgebra $GDer(A)$. 

\bigskip

Consider the universal Lie superalgebra $W(\Pi A)$ associated to the vector superspace $\Pi A$: this is the $\mathbb{Z}_+$-graded Lie superalgebra:
$$W(\Pi A) = \bigoplus_{k=-1}^{\infty} W_{k}(\Pi A)$$
where $W_{-1}=\Pi A$ and for all $k \geq 0$, $W_{k}(V) = \mathrm{Hom}(S^{k+1}(\Pi A),\Pi A)$ is the vector superspace of $(k+1)$-linear supersymmetric functions on $\Pi A$ with values in $\Pi A$. The Lie superalgebra structure on $W(\Pi A)$ is defined as follows: for $X \in W_p(\Pi A)$ and $Y \in W_{q}(\Pi A)$ with $p,q \geq -1$, we define $X \Box Y \in W_{p+q}(\Pi A)$ by:
\begin{equation}\label{boxproduct}
X \Box Y (a_0, \dots , a_{p+q}) = \sum_{{\tiny \begin{array}{c} i_0 < \dots < i_q \\ i_{q+1} < \dots < i_{q+p} \end{array}}} \epsilon_a(i_0, \dots , i_{p+q}) X(Y(a_{i_0}, \dots , a_{i_q}), a_{i_{q+1}}, \dots , a_{i_{q+p}}).
\end{equation}
Here $\epsilon_a(i_0, \dots , i_{p+q})= (-1)^{N}$ where $N$ is the number of interchanges of indices of odd $a_i$'s in the permutation $\sigma(s)=i_s$, $s= 0, \dots , p+q$. Then the bracket on $W(\Pi A)$ is given by:
$$[X,Y] = X \Box Y - (-1)^{\bar{p}(X) \bar{p}(Y)}Y \Box X.$$

As $GDer(A)$ is a subalgebra of the Lie superalgebra $W_0(\Pi A) = \mathrm{End}(\Pi A)$, we can consider its full prolongation ${\cal G}W^{as}(\Pi A)$: this is the $\mathbb{Z}_+$-graded subalgebra ${\cal G}W^{as}(\Pi A) = \bigoplus_{k = -1}^{\infty} {\cal G}W_{k}^{as}(\Pi A)$ of the Lie superalgebra $W(\Pi A)$ defined by setting ${\cal G}W_{-1}^{as}(\Pi A)=\Pi A$, ${\cal G}W_0^{as}(\Pi A) = GDer(\Pi A)$, and inductively for $k\geq 1$,
$${\cal G}W_{k}^{as}(\Pi A) = \{X \in W_{k}(\Pi A) \vert [X, W_{-1}(\Pi A)] \subset {\cal G}W_{k-1}^{as}(\Pi A)\}.$$

\begin{proposition}
For $k \geq 0$, the superspace ${\cal G}W_{k}^{as}(\Pi A)$ consists of linear maps $X: S^{k+1}(\Pi A) \rightarrow \Pi A$ satisfying the following generalized Leibniz rule:
\begin{equation}\label{genleibniz2}
X(a_0, \dots , a_{k-1}, bc) = X(a_0, \dots , a_{k-1}, b)c + (-1)^{p(b)p(c)}X(a_0, \dots , a_{k-1}, c)b - X(a_0, \dots , a_{k-1}, 1)bc
\end{equation}
for $a_0, \dots , a_{k-1}, b, c \in \Pi A$.
\end{proposition}
{\bf Proof.} According to formula (\ref{boxproduct}), for all $X \in W_{p}(\Pi A)$ and $Y \in W_{-1}(\Pi A) = \Pi A$, we have:
\begin{align}\label{equa3}
[X, Y](a_1, \dots , a_p) = X(Y,a_1, \dots , a_p)
\end{align}
with $a_1, \dots , a_p \in \Pi A$. Now we proceed by induction on $k \geq 0$: for $k=0$, ${\cal G}W_{0}^{as}(\Pi A) = GDer(A)$ and equality (\ref{genleibniz2}) holds by definition of  generalized derivation. Assume property (\ref{genleibniz2}) for elements in ${\cal G}W_{k-1}^{as}(\Pi A)$, and let $X$ in ${\cal G}W_{k}^{as}(\Pi A)$. For any $a_0, a_1, \dots , a_{k-1}, b, c \in \Pi A$, we have by (\ref{equa3}):
\begin{align*}
X(a_0, a_1, \dots , a_{k-1}, bc) = [X,a_0](a_1, \dots , a_{k-1}, bc).
\end{align*}
By definition of ${\cal G}W^{as}(\Pi A)$, we have $[X, a_0 ] \in {\cal G}W_{k-1}^{as}(\Pi A)$. Using the inductive hypothesis on $[X, a_0]$, we get:
\begin{align*}
[X,a_0](a_1, \dots , a_{k-1}, bc) =& [X,a_0](a_1, \dots , a_{k-1}, b)c + (-1)^{p(b)p(c)}[X,a_0](a_1, \dots , a_{k-1}, c)b\\ & - [X,a_0](a_1, \dots , a_{k-1}, 1)bc
\end{align*}
which is exactly formula (\ref{genleibniz2}) for $X$.
\hfill$\Box$

\bigskip

For $X \in \Pi W_{h-1}(\Pi A)$ and $Y \in \Pi W_{k-1}(\Pi A)$ with $h, k \geq 0$, we define their concatenation product $X \wedge Y \in \Pi W_{h+k-1}(\Pi A)$ by
\begin{equation}\label{concatproduct}
\begin{array}{ccc}
X \wedge Y(a_1, \dots a_{h+k}) & = & \displaystyle\sum_{{\tiny \begin{array}{c} i_1 < \dots < i_{h} \\ i_{h+1} < \dots < i_{h+k} \end{array}}} \epsilon_{a}(i_1, \dots , i_{h+k}) (-1)^{p(Y)(\bar{p}(a_{i_1})+ \dots + \bar{p}(a_{i_{h}}))}\\
&  & \,\,\,\,\,\,\,\,\,\,\,\,\,\,\,\,\,\,\,\,\,\,\,\,\,\,\,\times X(a_{i_1}, \dots , a_{i_{h}})Y(a_{i_{h+1}}, \dots , a_{i_{h+k}})
\end{array}
\end{equation}
where $\epsilon_a$ is defined as in (\ref{boxproduct}) with $a_1, \dots , a_{h+k} \in \Pi A$.

\begin{proposition}\label{UOGP}
$(\Pi {\cal G}W^{as}(\Pi A), [\cdot , \cdot ], \wedge)$ is a $\mathbb{Z}_+$-graded odd generalized Poisson superalgebra.
\end{proposition}

We will  denote $(\Pi {\cal G}W^{as}(\Pi A), [\cdot , \cdot ], \wedge)$ by ${\cal G}(A)$ and call it the universal odd generalized Poisson superalgebra associated to $A$. The rest of this section is devoted to the proof of Proposition \ref{UOGP}.

\begin{lemma}
$(\Pi {\cal G}W^{as}(\Pi A), \wedge)$ is a unital $\mathbb{Z}_+$-graded associative commutative superalgebra with parity $p$.
\end{lemma}
{\bf Proof.} It is already proved in \cite{DeSK} that $(\Pi W(\Pi A), \wedge)$ is a unital $\mathbb{Z}_+$-graded associative commutative superalgebra with parity $p$, therefore we only need to prove that for $X \in \Pi {\cal G}W^{as}_{h-1}(\Pi A)$ and $Y \in \Pi {\cal G}W_{k-1}^{as}(\Pi A)$ with $h, k \geq 0$, $X \wedge Y \in \Pi W_{h+k-1}(\Pi A)$ satisfies the generalized Leibniz rule (\ref{genleibniz2}). We have:
\begin{align}\label{equa1}
\begin{array}{c}
X \wedge Y(a_1, \dots , a_{h+k-1}, bc) = \\
\\
=\displaystyle\sum_{{\tiny \begin{array}{c} i_1 < \dots < i_{h} \\ i_{h+1} < \dots < i_{h+k}=h+k \end{array}}} \epsilon_{a_1, \dots , a_{h+k-1}, bc}(i_1, \dots , i_{h+k}) (-1)^{p(Y)(\bar{p}(a_{i_1})+ \dots + \bar{p}(a_{i_{h}}))}\\
\times X(a_{i_1}, \dots , a_{i_{h}})Y(a_{i_{h+1}}, \dots , a_{i_{h+k-1}}, bc)\\
\\
+ \displaystyle\sum_{{\tiny \begin{array}{c} i_1 < \dots < i_{h}=h+k \\ i_{h+1} < \dots < i_{h+k} \end{array}}} \epsilon_{a_1, \dots , a_{h+k-1}, bc}(i_1, \dots , i_{h+k}) (-1)^{p(Y)(\bar{p}(a_{i_1})+ \dots + \bar{p}(a_{i_{h-1}})+\bar{p}(bc))}\\
\times X(a_{i_1}, \dots , a_{i_{h-1}},bc)Y(a_{i_{h+1}}, \dots , a_{i_{h+k}})
\end{array}
\end{align}

\noindent For the first summand in the right hand side, since $i_{h+k} = h+k$, we have:
\begin{align*}
\epsilon_{a_1, \dots , a_{h+k-1}, bc}(i_1, \dots , i_{h+k}) &= \epsilon_{a_1, \dots , a_{h+k-1}, b}(i_1, \dots , i_{h+k})\\ &= \epsilon_{a_1, \dots , a_{h+k-1}, c}(i_1, \dots , i_{h+k})\\
&= \epsilon_{a_1, \dots , a_{h+k-1}, 1}(i_1, \dots , i_{h+k})
\end{align*}
and
\begin{align*}
Y(a_{i_{h+1}}, \dots , a_{i_{h+k-1}}, bc) =& Y(a_{i_{h+1}}, \dots , a_{i_{h+k-1}}, b)c + (-1)^{p(b)p(c)}Y(a_{i_{h+1}}, \dots , a_{i_{h+k-1}}, c)b\\ & - Y(a_{i_{h+1}}, \dots , a_{i_{h+k-1}}, 1)bc.
\end{align*}
In the second summand, since $i_{h+k} = h$, we have:
\begin{align*}
\epsilon_{a_1, \dots , a_{h+k-1}, bc}(i_1, \dots , i_{h+k}) &= \epsilon_{a_1, \dots , a_{h+k-1}, b}(i_1, \dots , i_{h+k}) (-1)^{p(c)(\bar{p}(a_{i_{h+1}})+ \dots + \bar{p}(a_{i_{h+k}}))} \\
& = \epsilon_{a_1, \dots , a_{h+k-1}, c}(i_1, \dots , i_{h+k})(-1)^{p(b)(\bar{p}(a_{i_{h+1}})+ \dots + \bar{p}(a_{i_{h+k}}))}\\
&= \epsilon_{a_1, \dots , a_{h+k-1}, 1}(i_1, \dots , i_{h+k})(-1)^{p(bc)(\bar{p}(a_{i_{h+1}})+ \dots + \bar{p}(a_{i_{h+k}}))}
\end{align*}
and
\begin{align*}
\begin{array}{l}
X(a_{i_1}, \dots , a_{i_{h-1}},bc)Y(a_{i_{h+1}}, \dots , a_{i_{h+k}}) =\\
\\ (-1)^{p(c)(p(Y)+\bar{p}(a_{i_{h+1}})+ \dots + \bar{p}(a_{i_{h+k}}))}X(a_{i_1}, \dots , a_{i_{h-1}},b)Y(a_{i_{h+1}}, \dots , a_{i_{h+k}})c\\
\\
+(-1)^{p(b)(p(c)+p(Y)+\bar{p}(a_{i_{h+1}})+ \dots + \bar{p}(a_{i_{h+k}}))}X(a_{i_1}, \dots , a_{i_{h-1}},c)Y(a_{i_{h+1}}, \dots , a_{i_{h+k}})b\\
\\
- (-1)^{p(bc)(p(Y)+\bar{p}(a_{i_{h+1}})+ \dots + \bar{p}(a_{i_{h+k}}))}X(a_{i_1}, \dots , a_{i_{h-1}},1)Y(a_{i_{h+1}}, \dots , a_{i_{h+k}})bc
\end{array}
\end{align*}
The generalized Leibniz rule for $X \wedge Y$ then follows by replacing these equalities in (\ref{equa1}).
\hfill$\Box$

\bigskip

\noindent It remains to prove that the Lie bracket on $\Pi {\cal G}W^{as}(\Pi A)$ satisfies the generalized odd Leibniz rule (\ref{oddleibniz}). This follows from the following lemma.

\begin{lemma}
The following equalities hold for $X,Y,Z\in \Pi {\cal G}W^{as}(\Pi A)$:
\begin{align*}
X \Box ( Y \wedge Z) &= (X \Box Y) \wedge Z + (-1)^{\bar{p}(X)p(Y)} Y \wedge (X \Box Z) - (X \Box 1)\wedge Y \wedge Z,\\
(X \wedge Y) \Box Z &= X \wedge (Y \Box Z) + (-1)^{p(Y)\bar{p}(Z)}(X \Box Z) \wedge Y.
\end{align*}
\end{lemma}
{\bf Proof.} An analogue result is proved in \cite[Lemma 3.5]{DeSK}.  For $X \in \Pi {\cal G}W_{l-k}^{as}(\Pi A)$, $Y \in \Pi {\cal G}W_{h-1}^{as}(\Pi A)$ and $Z \in \Pi {\cal G}W_{k-h-1}^{as}(\Pi A)$ with $h, k-h, l-k+1\geq 0$, we have:
\begin{align}\label{equa2}
\begin{array}{rl}
X \Box ( Y \wedge Z)(a_{1}, \dots , a_{l}) =& \displaystyle \sum_{\begin{tiny}\begin{array}{c} i_1 < \dots < i_{h} \\ i_{h+1} < \dots < i_{k} \\ i_{k+1} < \dots < i_{l} \end{array}
\end{tiny}} \epsilon_a(i_1, \dots , i_l)(-1)^{p(Z)(\bar{p}(a_{i_{1}})+ \dots + \bar{p}(a_{i_{h}}))} \\
& \times X(Y(a_{i_1}, \dots , a_{i_{h}})Z(a_{i_{h+1}}, \dots , a_{i_{k}}), a_{i_{k+1}}, \dots , a_{i_{l}})
\end{array}
\end{align}
The generalized Leibniz rule for $X$ can be rewritten in the following way:
\begin{align*}
X(bc, a_{k+1}, \dots , a_{l}) =& (-1)^{p(c)(\bar{p}(a_{k+1}) + \dots + \bar{p}(a_{l}))}X(b, a_{k+1}, \dots , a_{l})c \\
&+ (-1)^{p(b)\bar{p}(X)} b X(c, a_{k+1}, \dots , a_{l})\\
&-(-1)^{p(bc)(\bar{p}(a_{k+1}) + \dots + \bar{p}(a_{l}))}X(1, a_{k+1}, \dots , a_{l})bc.
\end{align*}
Using this equality in (\ref{equa2}), $X \Box ( Y \wedge Z)(a_{1}, \dots , a_{l})$ is then of the form:
$$X \Box ( Y \wedge Z)(a_{1}, \dots , a_{l}) = A + B - C.$$
The first term $A$ is equal to
\begin{align*}
\begin{array}{c}
\displaystyle \sum_{\begin{tiny}\begin{array}{c} i_1 < \dots < i_{h} \\ i_{h+1} < \dots < i_{k} \\ i_{k+1} < \dots < i_{l} \end{array}
\end{tiny}} \epsilon_a(i_1, \dots , i_l)(-1)^{p(Z)(\bar{p}(a_{i_{1}})+ \dots + \bar{p}(a_{i_{h}}))} (-1)^{(p(Z)+ \bar{p}(a_{i_{h+1}})+ \dots + \bar{p}(a_{i_{k}}))(\bar{p}(a_{i_{k+1}})+ \dots + \bar{p}(a_{i_{l}}))} \\
 \times X(Y(a_{i_1}, \dots , a_{i_{h}}), a_{i_{k+1}}, \dots , a_{i_{l}})Z(a_{i_{h+1}}, \dots , a_{i_{k}})=
\end{array}
\end{align*}
\begin{align*}
\begin{array}{c}
=\displaystyle \sum_{\begin{tiny}\begin{array}{c} i_1 < \dots < i_{h} \\ i_{h+1} < \dots < i_{l-k+h} \\ i_{l-k+h+1} < \dots < i_{l} \end{array}
\end{tiny}} \epsilon_a(i_1, \dots , i_l)(-1)^{(\bar{p}(a_{i_{h+1}}) + \dots + \bar{p}(a_{i_{l-k+h}}))(\bar{p}(a_{i_{l-k+h+1}}) + \dots + \bar{p}(a_{i_{l}}))}\\
\times (-1)^{p(Z)(\bar{p}(a_{i_{1}})+ \dots + \bar{p}(a_{i_{h}}))} (-1)^{(p(Z)+ \bar{p}(a_{i_{l-k+h+1}})+ \dots + \bar{p}(a_{i_l}))(\bar{p}(a_{i_{h+1}})+ \dots + \bar{p}(a_{i_{l-k+h}}))} \\
\\ \times X(Y(a_{i_1}, \dots , a_{i_{h}}), a_{i_{h+1}}, \dots , a_{i_{l-k+h}})Z(a_{i_{l-k+h+1}}, \dots , a_{i_{l}})=
\end{array}
\end{align*}
$$=(X \Box Y) \wedge Z(a_{1}, \dots , a_{l}).$$

\medskip

The second term $B$ is equal to
\begin{align*}
\begin{array}{c}
\displaystyle \sum_{\begin{tiny}\begin{array}{c} i_1 < \dots < i_{h} \\ i_{h+1} < \dots < i_{k} \\ i_{k+1} < \dots < i_{l} \end{array}
\end{tiny}} \epsilon_a(i_1, \dots , i_l)(-1)^{p(Z)(\bar{p}(a_{i_{1}})+ \dots + \bar{p}(a_{i_{h}}))} (-1)^{(p(Y)+ \bar{p}(a_{i_{1}})+ \dots + \bar{p}(a_{i_{h}}))\bar{p}(X)} \\
 \times Y(a_{i_1}, \dots , a_{i_{h}})X(Z(a_{i_{h+1}}, \dots , a_{i_{k}}), a_{i_{k+1}}, \dots , a_{i_{l}})\\
\\
= (-1)^{p(Y)\bar{p}(X)} \displaystyle \sum_{\begin{tiny}\begin{array}{c} i_1 < \dots < i_{h} \\ i_{h+1} < \dots < i_{k} \\ i_{k+1} < \dots < i_{l} \end{array}
\end{tiny}} \epsilon_a(i_1, \dots , i_l)(-1)^{(p(Z)+\bar{p}(X))(\bar{p}(a_{i_{1}})+ \dots + \bar{p}(a_{i_{h}}))} \\
\\
 \times Y(a_{i_1}, \dots , a_{i_{h}})X(Z(a_{i_{h+1}}, \dots , a_{i_{k}}), a_{i_{k+1}}, \dots , a_{i_{l}})=(-1)^{\bar{p}(X)p(Y)} Y \wedge (X \Box Z)(a_{1}, \dots , a_{l})
\end{array}
\end{align*}
since $p(X \Box Z) = \bar{p}(X) + p(Z)$.

\medskip

Finally, the third term $C$ is equal to
\begin{align*}
\begin{array}{c}
\displaystyle \sum_{\begin{tiny}\begin{array}{c} i_1 < \dots < i_{h} \\ i_{h+1} < \dots < i_{k} \\ i_{k+1} < \dots < i_{l} \end{array}
\end{tiny}} \epsilon_a(i_1, \dots , i_l)(-1)^{p(Z)(\bar{p}(a_{i_{1}})+ \dots + \bar{p}(a_{i_{h}}))} (-1)^{(p(Y)+p(Z)+ \bar{p}(a_{i_1})+ \dots + \bar{p}(a_{i_{k}}))(\bar{p}(a_{i_{k+1}})+ \dots + \bar{p}(a_{i_l}))} \\
 \times X(1, a_{i_{k+1}}, \dots , a_{i_{l}}) Y(a_{i_1}, \dots , a_{i_{h}})Z(a_{i_{h+1}}, \dots , a_{i_{k}})\\
\\
=\displaystyle \sum_{\begin{tiny}\begin{array}{c} i_1 < \dots < i_{h} \\ i_{h+1} < \dots < i_{k} \\ i_{k+1} < \dots < i_{l} \end{array}
\end{tiny}} \epsilon_a(i_1, \dots , i_l)(-1)^{(\bar{p}(a_{i_1})+ \dots + \bar{p}(a_{i_{k}}))(\bar{p}(a_{i_{k+1}})+ \dots + \bar{p}(a_{i_l}))}\\
\times (-1)^{p(Z)(\bar{p}(a_{i_{1}})+ \dots + \bar{p}(a_{i_{h}}) + \bar{p}(a_{i_{k+1}})+ \dots + \bar{p}(a_{i_l}))} (-1)^{p(Y)(\bar{p}(a_{i_{k+1}})+ \dots + \bar{p}(a_{i_l}))} \\
\\
 \times X(1, a_{i_{k+1}}, \dots , a_{i_{l}}) Y(a_{i_1}, \dots , a_{i_{h}})Z(a_{i_{h+1}}, \dots , a_{i_{k}})=
\end{array}
\end{align*}

\begin{align*}
=\displaystyle \sum_{\begin{tiny}\begin{array}{c} i_1 < \dots < i_{l-k} \\ i_{l-k+1} < \dots < i_{l-k+h} \\ i_{l-k+h+1} < \dots < i_{l} \end{array}
\end{tiny}} \epsilon_a(i_1, \dots , i_l)
(-1)^{p(Z)(\bar{p}(a_{i_{1}})+ \dots + \bar{p}(a_{i_{l-k+h}}))} (-1)^{p(Y)(\bar{p}(a_{i_{1}})+ \dots + \bar{p}(a_{i_{l-k}}))} \\
\\
 \times (X \Box 1)(a_{i_{1}}, \dots , a_{i_{l-k}}) Y(a_{i_{l-k+1}}, \dots , a_{i_{l-k+h}})Z(a_{i_{l-k+h+1}}, \dots , a_{i_{l}})=
(X \Box 1) \wedge Y \wedge Z (a_1, \dots , a_l).
\end{align*}
 This proves the first equality. The second equality can be proved in the same way, using the definition of the box product (\ref{boxproduct}) and the concatenation product (\ref{concatproduct}).
\hfill$\Box$

\section{The main construction}
Let $({\cal N},\{\cdot,\dots,\cdot\},\cdot)$ be a generalized $n$-Nambu-Poisson algebra and denote by $\Pi {\cal N}$ the space ${\cal N}$ with reversed parity. Define
\begin{equation}
\begin{array}{c}
\mu: \Pi {\cal N}\otimes\dots \otimes \Pi {\cal N} \rightarrow \Pi {\cal N}\\
\mu(f_1, \dots, f_n)=\{f_1, \dots, f_n\}.
\end{array}
\label{mu}
\end{equation}
Then $\mu$ is a supersymmetric function on $(\Pi {\cal N})^{\otimes n}$ \cite[Lemma 1.2]{CantaK3}. 
Furthermore $\mu$ satisfies the generalized Leibniz rule 
$$\mu(f_1, \dots, f_{n-1},gh)= \mu(f_1, \dots, f_{n-1},g)h+g\mu(f_1, \dots, f_{n-1},h)-\mu(f_1,\dots, f_{n-1},1)gh,$$
hence $\mu$ lies in ${\cal G}W^{as}_{n-1}(\Pi {\cal N})$. 

Let $OP({\cal N})$ be the odd Poisson subalgebra of ${\cal G}({\cal N})$ generated by $\Pi {\cal N}$ and $\mu$. Then, by construction,
$OP({\cal N})$ is a transitive Lie subalgebra of ${\cal G}W^{as}(\Pi {\cal N})$, hence it is a transitive subalgebra of $W(\Pi {\cal N})$. Furthermore
$OP({\cal N})$ is a $\Z_+$-graded odd Poisson subalgebra of  ${\cal G}({\cal N})$. Let us denote by $OP({\cal N})=\oplus_{j\geq -1}{\mathcal P}_j({\cal N})$
its depth 1 $\Z$-grading as a Lie superalgebra.
\begin{proposition}\label{simple}
If ${\cal N}$ is a simple generalized $n$-Nambu-Poisson algebra then $OP({\cal N})$ is a simple generalized odd Poisson superalgebra.
\end{proposition}
{\bf Proof.} Let $I$ be a non-zero ideal of $OP({\cal N})$. Then, by transitivity, $I\cap {\mathcal P}_{-1}({\cal N})=I\cap {\cal N}\neq 0$.
Note that $I\cap {\cal N}$ is a Nambu-Poisson ideal of ${\cal N}$. Indeed, $(I\cap {\cal N})\cdot N=(I\cap {\cal N})\wedge {\cal N}\subset I\cap {\cal N}$ and $[I\cap {\cal N}, {\cal N}]
\subset [{\cal N},{\cal N}]=0$. Since ${\cal N}$ is simple, $I\cap {\cal N}={\cal N}$, hence $1\in I$, hence $I=OP({\cal N})$. \hfill$\Box$

\begin{remark}\em
We recall that since $({\cal N}, \{\cdot,\dots,\cdot\})$ is an $n$-Lie algebra, the Filippov-Jacobi identity holds, i.e., for every $a_1, \dots, a_{n-1}
\in {\cal N}$, the map $D_{a_1, \dots, a_{n-1}}:{\cal N}\rightarrow {\cal N}$, $D_{a_1, \dots, a_{n-1}}(a)=\{a_1,\dots, a_{n-1},a\}$
is a derivation of $({\cal N}, \{\cdot,\dots,\cdot\})$. By \cite[Lemma 2.1(b)]{CantaK4}, this is equivalent to the condition $[\mu, D_{a_1, \dots, a_{n-1}}]=0$ in $OP({\cal N})$.
By (\ref{mu}), we have: $D_{a_1, \dots, a_{n-1}}=[[\mu, a_1], \dots, a_{n-1}]$, therefore $\mu$ satisfies the following condition:
$$[\mu, [[\mu, a_1], \dots, a_{n-1}]]=0 \,\,\,\, {\mbox{for every}}\, a_1, \dots, a_{n-1}
\in {\cal N}.$$
\end{remark}

\begin{definition}\label{GoodPair}
We say that a pair $({\cal P},\mu)$, consisting of a $\Z_+$-graded 
generalized odd Poisson superalgebra ${\cal P}$ and an element $\mu\in {\mathcal P}_{n-1}$ of parity $p(\mu)\equiv n$ (mod $2$), is a good $n$-pair if it satisfies 
the following properties:
\begin{itemize}
\item[G1)] ${\cal P}=\oplus_{j\geq -1}{\mathcal P}_j$ is a transitive $\Z$-graded Lie superalgebra of depth 1 such that ${\cal P}_{-1}$ is completely odd;
\item[G2)] $\mu$ and ${\cal P}_{-1}$ generate ${\cal P}$ as a (generalized) odd Poisson superalgebra;
\item[G3)] $[\mu, [[\mu, a_1], \dots, a_{n-1}]]=0$ for every $a_1, \dots, a_{n-1}
\in {\cal P}_{-1}$.
\end{itemize}
\end{definition} 

\begin{example}\label{PO2}\em
Let ${\cal P}=PO(2h,2h)$, $h\geq 1$, with the grading of type $(0,\dots,0|1, \dots,1)$, and let $\mu=\sum_{i=1}^h\xi_i\xi_{i+h}$. Then 
$({\cal P}, \mu)$ is a good $2$-pair. Indeed, for $1\leq i\leq h$, $[x_i, \mu]_{HO}=\xi_{h+i}$ and $[x_{h+i}, \mu]_{HO}=-\xi_{i}$,
therefore ${\cal P}_{-1}$ and $\mu$ generate ${\cal P}$.
Furthermore, for $f\in{\cal P}_{-1}=\F[[x_1, \dots, x_n]]$, we have:
$[f,\mu]_{HO}=\sum_{i=1}^h(\frac{\partial f}{\partial x_i}\xi_{i+h}-\frac{\partial f}{\partial x_{i+h}}\xi_{i})$, hence
$$[\mu, [f,\mu]_{HO}]_{HO}=\sum_{i,j=1}^h[\xi_j\xi_{j+h}, \frac{\partial f}{\partial x_i}\xi_{i+h}-\frac{\partial f}{\partial x_{i+h}}\xi_{i}]_{HO}=$$
$$=\sum_{i,j=1}^h(\frac{\partial^2 f}{\partial x_i\partial x_j}\xi_{j+h}\xi_{i+h}-\frac{\partial^2 f}{\partial x_{j+h}\partial x_{i}}\xi_j\xi_{i+h}
-\frac{\partial^2 f}{\partial x_j\partial x_{i+h}}\xi_{j+h}\xi_{i}+\frac{\partial^2 f}{\partial x_{i+h}\partial x_{j+h}}\xi_j\xi_{i})=0.$$
Therefore $({\cal P}, \mu)$ satisfies property $G3)$.
\end{example}

\begin{example}\label{POn}\em
Let ${\cal P}=PO(n,n)$ with the grading of type $(0,\dots,0|1, \dots,1)$, and let $\mu=\xi_1\dots\xi_n$. Then 
$({\cal P}, \mu)$ is a good $n$-pair. Indeed, $[x_{n-1},[\dots, [x_2, [x_1, \mu]]]]_{HO}=\xi_n$,
and, similarly all the $\xi_i$'s can be obtained by commuting $\mu$ with different $x_j$'s. 
Therefore ${\cal P}_{-1}$ and $\mu$ generate ${\cal P}$.
Furthermore, let $f=\sum_{i=1}^nf_i\xi_i\in{\cal P}_0$, with $f_i\in\F[[x_1, \dots, x_n]]$, such that
\begin{equation}
\sum_{i=1}^n\frac{\partial f_i}{\partial x_i}=0.
\label{0div}
\end{equation}
Then
$[f,\mu]_{HO}=\sum_{i=1}^n\frac{\partial f_i}{\partial x_i}\xi_1\dots\xi_n=0.$
Notice that all elements of the form $[[\mu, a_1], \dots, a_{n-1}]$ with $a_1, \dots, a_{n-1}
\in {\cal P}_{-1}=\F[[x_1, \dots, x_n]]$ satisfy property (\ref{0div}), hence
$({\cal P}, \mu)$ satisfies property $G3)$.
\end{example}

\begin{example}\label{PO2h+1}\em Let ${\cal P}=PO(2h+1, 2h+2)$, $h\geq 1$, with the grading of type $(0, \dots, 0|1, \dots, 1,1)$, and let
$\mu=\sum_{i=1}^{h+1}\xi_i\xi_{i+h+1}$ (recall that $\xi_{2h+2}=\tau$). Then $({\cal P}, \mu)$ is a good $2$-pair. Indeed, we have:
$[1, \mu]_{KO}=2\xi_{h+1}$ and $[x_i, \mu]_{KO}=\xi_{i+h+1}-x_i\xi_{h+1}$ for $1\leq i\leq h+1$, $[x_{i+h+1}, \mu]_{KO}=-\xi_{i}-x_{i+h+1}\xi_{h+1}$
for $1\leq i\leq h$.
Hence ${\cal P}_{-1}$ and $\mu$ generate ${\cal P}$.
Furthermore, if $f\in {\cal P}_{-1}=\F[[x_1, \dots, x_n]]$, we have:
$$[f,\mu]_{KO}=\sum_{i=1}^h(\frac{\partial f}{\partial x_i}\xi_{i+h+1}-\frac{\partial f}{\partial x_{i+h+1}}\xi_{i})+
\frac{\partial f}{\partial x_{h+1}}\xi_{2h+2}-(E-2)(f)\xi_{h+1},$$ hence
$$[\mu,[f,\mu]_{KO}]_{KO}=[\sum_{j=1}^{h+1}\xi_j\xi_{h+1+j}, \sum_{i=1}^h(\frac{\partial f}{\partial x_i}\xi_{i+h+1}-\frac{\partial f}{\partial x_{i+h+1}}\xi_{i})+
\frac{\partial f}{\partial x_{h+1}}\xi_{2h+2}-(E-2)(f)\xi_{h+1}]_{KO}$$
$$=\sum_{
i=1, \dots, h;\\ j=1, \dots, h+1
}\xi_{h+1+j}(\frac{\partial^2 f}{\partial x_j\partial x_i}\xi_{h+i+1}-
\frac{\partial^2 f}{\partial x_j\partial x_{i+h+1}}\xi_{i})+\sum_{j=1}^{h+1}\xi_{h+1+j}
(\frac{\partial^2 f}{\partial x_j\partial x_{h+1}}\xi_{2h+2}$$
$$-\frac{\partial ((E-2)(f)}{\partial x_j}\xi_{h+1})
-\sum_{i,j=1}^h\xi_j(\frac{\partial^2 f}{\partial x_{h+1+j}\partial x_i}\xi_{h+i+1}-
\frac{\partial^2 f}{\partial x_{h+1+j}\partial x_{i+h+1}}\xi_{i})-\sum_{j=1}^{h}\xi_{j}
(\frac{\partial^2 f}{\partial x_{h+1+j}\partial x_{h+1}}\xi_{2h+2}$$
$$-\frac{\partial ((E-2)(f)}{\partial x_{h+1+j}}\xi_{h+1})-\xi_{h+1}(E-2)([f,\mu]_{KO})
=\xi_{2h+2}\sum_{i=1}^h(\frac{\partial^2f}{\partial x_{h+1}\partial x_i}\xi_{h+i+1}-
\frac{\partial^2f}{\partial x_{h+1}\partial x_{i+h+1}}\xi_{i})$$
$$+\sum_{j=1}^{h+1}\frac{\partial^2f}{\partial x_{j}\partial x_{h+1}}\xi_{h+1+j}\xi_{2h+2}-
\sum_{j=1}^{h+1}\xi_{h+1+j}\frac{\partial ((E-2)(f))}{\partial x_{j}}\xi_{h+1}-
\sum_{j=1}^h(\frac{\partial^2f}{\partial x_{h+j+1}\partial x_{h+1}}\xi_j\xi_{2h+2}$$
$$-\frac{\partial ((E-2)(f))}{\partial x_{h+1+j}}\xi_{j}\xi_{h+1})
-\xi_{h+1}(E-2)([f,\mu]_{KO})=\xi_{h+1}(\sum_{j=1}^{h+1}\xi_{h+1+j}\frac{\partial((E-2)(f))}{\partial x_j}$$
$$-\sum_{j=1}^{h}\frac{\partial((E-2)(f))}{\partial x_{h+1+j}}\xi_j
-(E-2)(\sum_{i=1}^h(\frac{\partial f}{\partial x_i}\xi_{i+h+1}-\frac{\partial f}{\partial x_{i+h+1}}\xi_{i})+
\frac{\partial f}{\partial x_{h+1}}\xi_{2h+2}))
=0.$$
\end{example}

\begin{example}\label{POntau}\em Let ${\cal P}=PO(n,n+1)$, with the grading of type $(0, \dots, 0|1, \dots, 1,1)$, and let
$\mu=\xi_1\dots\xi_n\tau$ (recall that $\tau=\xi_{n+1}$). Then $({\cal P}, \mu)$ is a good $(n+1)$-pair. Indeed, we have:
$[1, \mu]_{KO}=2(-1)^{n+1}\xi_1\dots\xi_n$, $[x_{i_1},[\dots,[x_{i_{n-1}}, \xi_1\dots\xi_n]_{KO}]_{KO}]_{KO}=\pm \xi_{i_n}$ for $i_1\neq \dots\neq i_{n-1}\neq i_n$,
 $[x_i, \xi_i\dots\xi_n\tau]_{KO}=\xi_{i+1}\dots\xi_n\tau+(-1)^{n-i}x_i\xi_i\dots\xi_n$ for $1\leq i\leq n$.
Hence ${\cal P}_{-1}$ and $\mu$ generate ${\cal P}$.
Now let $div_1=\Delta+(E-n)\frac{\partial}{\partial \tau}$ where $\Delta=\sum_{i=1}^n
\frac{{\partial}^2}{\partial x_i\partial\xi_i}$ is the odd Laplacian, and let
$f=\sum_{i=1}^{n+1}f_i\xi_i\in {\cal P}_0$, $f_i\in\F[[x_1, \dots,x_n]]$, such that
$0=div_1(f)=\sum_{i,j=1}^n\frac{\partial f_i}{\partial x_j}+(E-n)(f_{n+1})$.
Then we have:
$$[\sum_{i=1}^{n+1}f_i\xi_i, \mu]_{KO}=[\sum_{i=1}^{n+1}f_i\xi_i, \mu]_{HO}+
\sum_{i=1}^{n+1}(E-2)(f_i\xi_i)(-1)^n\xi_1\dots\xi_n-f_{n+1}(n-2)\xi_1\dots\xi_n\tau$$
$$=\sum_{i,j=1}^n\frac{\partial f_i}{\partial x_j}\xi_1\dots\xi_n\tau+(-1)^n(E-2)(f_{n+1}\xi_{n+1})\xi_1\dots\xi_n
-(n-2)f_{n+1}\xi_1\dots\xi_n\tau$$
$$=(\sum_{i,j=1}^n\frac{\partial f_i}{\partial x_j}+(E-2)(f_{n+1})-(n-2)f_{n+1})\xi_1\dots\xi_n\tau=0.$$
Notice that, since $div_1(\mu)=0$ and $div_1(f)=0$ for every $f\in {\cal P}_{-1}$, then $div_1([[[\mu, a_1], \dots, a_{n-1}]])=0$ for every $a_1, \dots, a_{n-1}
\in {\cal P}_{-1}$. Hence property $G3)$ is satisfied.
\end{example}

\begin{remark}\label{GaugeGood}\em Let us consider ${\cal P}=PO(k,k)$ (resp.\ ${\cal P}=PO(k,k+1)$) with the grading of type 
$(0,\dots,0|1,\dots,1)$ (resp.\ $(0,\dots,0|1,\dots,1,1)$). Let $\varphi\in {\cal P}_{-1}$ be an invertible element. By Remark \ref{gaugegradings}, the grading of type $(0,\dots,0|1,\dots,1)$  (resp.\ $(0,\dots,0|1,\dots,1,1)$) 
defines a $\Z_+$-graded structure on the odd generalized Poisson
superalgebra ${\cal P}^{\varphi}$, such that ${\cal P}_j={\cal P}_j^{\varphi}$. 
Then, by (\ref{gaugequiv}),  $({\cal P}, \mu)$ is a good $n$-pair with respect to this grading if and only if $({\cal P}^{\varphi}, \varphi^{-1}\mu)$ is. 
\end{remark}

\medskip

The map ${\cal N} \mapsto (OP({\cal N}), \mu)$ establishes a correspondence between 
(simple) generalized $n$-Nambu-Poisson algebras ${\cal N}$ and good $n$-pairs $(OP({\cal N}), \mu)$. 
We now want to show that this correspondence  is bijective.
\begin{lemma}\label{0component} Let ${\cal N}$ be a generalized $n$-Nambu-Poisson algebra. Then 
the $0$-th graded component ${\mathcal P}_0(N)$ of $OP(N)$ is generated, as a Lie
superalgebra, by elements of the form
$$[a_1, [a_2,\dots, [a_{n-1}, \mu]]]b$$
with $a_i,b\in\Pi {\cal N}$.
\end{lemma}
{\bf Proof.} Let $L_{-1}:=\Pi {\cal N}$ and let $L_0$ be the
Lie subsuperalgebra of ${\cal G}W^{as}_0(\Pi {\cal N})=GDer(\Pi {\cal N})$ generated by the elements 
of the form $[a_1, [a_2, \dots, [a_{n-1}, \mu b]]]$ with $a_1, \dots, a_{n-1}, b\in \Pi {\cal N}$.
Note that, since  ${\cal G}W^{as}_0(\Pi {\cal N})$ is $\Z$-graded of depth 1, and $1\in {\cal N}$, the restriction to ${\cal N}$ of the derivation $D$
of ${\cal G}({\cal N})$  is zero, hence
$$[a_1, [a_2, \dots, [a_{n-1}, \mu b]]]=[a_1, [a_2, \dots, [a_{n-1}, \mu ]]]b.$$

An induction argument on the length of the commutators of the generating elements of $L_0$ shows
 that $L_0$ is stable with respect to the concatenation product by elements of $\Pi{\cal N}$. 

Let $L$ be the full prolongation of $L_{-1}\oplus L_0$, i.e., $L=L_{-1}\oplus L_0\oplus (\oplus_{j\geq 1}L_j)$,
where $L_j=\{\varphi\in {\cal G}W^{as}(\Pi {\cal N})~|~[\varphi, L_{-1}]\subset L_{j-1}\}$. Note that $L_j$, for $j\geq 1$, is stable with respect to
the concatenation product by elements of $\Pi{\cal N}$. Indeed, if $\varphi\in L_j$, then
$$[\varphi \Pi{\cal N}, L_{-1}]=[\varphi, L_{-1}]\Pi{\cal N}\subset L_{j-1}\Pi{\cal N},$$
hence one can conclude by induction on $j$ since $L_0\Pi{\cal N}\subset L_0$. It follows that $L$ is closed under the concatenation
product, hence it is an odd generalized Poisson subsuperalgebra of ${\cal G}W^{as}(\Pi {\cal N})$. Indeed, using induction on $i+j\geq 0$, one shows that
$L_iL_j\subset L$ for every $i,j\geq 0$.

It follows that $OP({\cal N})$ is an odd generalized subsuperalgebra of $L$, since $L$ is an odd generalized Poisson superalgebra containing $\Pi {\cal N}$ and
$\mu$. As a consequence, the  $0$-th graded component ${\mathcal P}_0({\cal N})$ of $OP({\cal N})$ is generated, as a Lie
superalgebra, by elements of the form
$$[a_1, [a_2,\dots, [a_{n-1}, \mu]]]b$$
with $a_i,b\in\Pi {\cal N}$.
\hfill$\Box$

\begin{proposition}\label{reversedarrow}
Let $({\cal P}, \mu)$ be a good $n$-pair, and define on   ${\cal N}:=\Pi {\cal P}_{-1}$ the following product:
$$\{x_1, \dots, x_n\}=[\dots [[\mu, x_1], \dots, x_n]].$$ Then:
\begin{itemize}
\item[(a)] $({\cal N}, \{\cdot, \dots, \cdot\}, \wedge)$ is a generalized Nambu-Poisson algebra, $\wedge$ being the restriction  to ${\cal N}$ of the 
commutative associative product $\wedge$ defined on ${\cal P}$.
\item[(b)] If ${\cal P}$ is a simple odd generalized Poisson superalgebra, then $({\cal N}, \{\cdot, \dots, \cdot\}, \wedge)$
is a simple generalized Nambu-Poisson algebra.
\end{itemize}
\end{proposition}
{\bf Proof.} $(a)$ By Definitions \ref{OGP} and \ref{graded}, ${\cal N}={\cal Q}_0$ is a commutative associative subalgebra of ${\cal P}$.
Furthermore $\{\cdot, \dots, \cdot\}$ is an $n$-Lie bracket due to \cite[Prop. 2.4]{CantaK4} and property $G3)$. Finally,
for $f_1, \dots, f_{n-1}, g, h\in \Pi{\cal P}_{-1}$, we have:

\noindent
$\{f_1, \dots, f_{n-1}, gh\}=[[\dots [\mu, f_1], \dots, f_{n-1}], gh]=
[[\dots [\mu, f_1], \dots, f_{n-1}], g]h+g[[\dots [\mu, f_1], \dots, f_{n-1}], h]$

\medskip
\noindent
$+(-1)^{p([\dots [\mu, f_1], \dots, f_{n-1}])+1}[1,[\dots [\mu, f_1], \dots, f_{n-1}]]gh=
\{f_1, \dots, f_{n-1}, g\}h+g\{f_1, \dots, f_{n-1}, h\}$

\medskip
\noindent
$-(-1)^{p([\dots [\mu, f_1], \dots, f_{n-1}])+1}
(-1)^{\bar{p}([\dots [\mu, f_1], \dots, f_{n-1}])}\{f_1, \dots, f_{n-1},1\}gh$

\medskip
\noindent
$=\{f_1, \dots, f_{n-1}, g\}h+g\{f_1, \dots, f_{n-1}, h\}-\{f_1, \dots, f_{n-1},1\}gh$.

$(b)$ Now we want to show that if ${\cal P}$ is simple, then ${\cal N}$ is simple. Suppose that $I$ is a non zero ideal of ${\cal N}$,
and let $\tilde{I}$ be the ideal of ${\cal P}$ generated by $\Pi I$ and $\mu$: $\tilde{I}=\oplus_{j\geq -1}\tilde{I}_j$,
with $\tilde{I}_j\subset {\cal P}_j$. We want to show that $\tilde{I}_{-1}=\tilde{I}\cap {\cal P}_{-1}=\Pi I$.
In fact, the concatenation product by elements in $\oplus_{j\geq 1}{\cal Q}_j$ maps ${\cal Q}_0$ to $\oplus_{j\geq 1}{\cal Q}_j$
hence it does not produce any element in ${\cal P}_{-1}={\cal Q}_0$.
On the other hand, $I\wedge {\cal Q}_0=I\wedge {\cal N}\subset I$ since $I$ is an ideal of ${\cal N}$. The bracket between elements in 
$\oplus_{j\geq 0}{\cal P}_j$ lies in $\oplus_{j\geq 0}{\cal P}_j$ and the bracket between $I$ and elements in $\oplus_{j\geq 1}{\cal P}_j$
lies in $\oplus_{j\geq 0}{\cal P}_j$. Therefore we just need to consider the brackets between elements in $I$ and elements in
${\cal P}_0$. By hypothesis, ${\cal P}$ is generated by ${\cal P}_{-1}$ and $\mu$, hence, by the same argument as in Lemma \ref{0component},
${\cal P}_0$ is generated by elements of the form $[a_1,[a_2,\dots, [a_{n-1}, \mu]]]b$ with $a_i, b\in \Pi {\cal P}_{-1}$. We have:
$$[I, [a_1,[a_2,\dots, [a_{n-1}, \mu]]]b]= [I, [a_1,[a_2,\dots, [a_{n-1}, \mu]]]]b$$
since $[I,b]=0$ and $D|_I=0$. Since $[I, [a_1,[a_2,\dots, [a_{n-1}, \mu]]]]=\{I, a_1, \dots, a_{n-1}\}$ and $I$ is an ideal of
${\cal N}$, $[I, {\cal P}_0]\subset I$. \hfill$\Box$

\medskip

\begin{definition} Two good $n$-pairs $({\cal P}, \mu)$ and $({\cal P}', \mu')$ are called isomorphic if there exists an odd
Poisson superalgebras isomorphism $\Phi: {\cal P} \rightarrow {\cal P}'$ such that $\Phi({\cal P}_j)={\cal P}'_j$,
$\Phi({\cal Q}_j)={\cal Q}'_j$ for all $j$ and $\phi(\mu)\in\F^{\times}\mu'$.
\end{definition}

\begin{theorem}\label{bijection} The map
$${\cal N}\rightarrow (OP({\cal N}), \mu)$$
with $\mu$ defined as in (\ref{mu}), establishes a bijection between isomorphism classes of generalized $n$-Nambu-Poisson algebras and
isomorphism classes of good $n$-pairs. Moreover:
\begin{itemize}
\item[(i)] ${\cal N}$ is simple (linearly compact) if and only if $OP({\cal N})$ is;
\item[(ii)] ${\cal N}$ is a Nambu-Poisson algebra if and only if $OP({\cal N})$ is an odd Poisson superalgebra.
\end{itemize}
\end{theorem}
{\bf Proof.} The proof follows immediately from Propositions \ref{simple} and \ref{reversedarrow}. The fact that the linear compactness
of ${\cal N}$ implies that of $OP({\cal N})$ can be proved in the same way as in \cite[Proposition 2.4]{CantaK4}.
\hfill$\Box$

\begin{remark}\label{correspondence}\em One can check (see also \cite{CantaK4}) that if ${\cal N}$ is the
$n$-Nambu algebra, then $(OP({\cal N}), \mu)=(PO(n,n),\xi_1\dots\xi_n)$ and if ${\cal N}$ is the
$n$-Dzhumaldidaev algebra, then $(OP({\cal N}), \mu)=(PO(n-1,n),\xi_1\dots\xi_{n-1}\tau)$.
\end{remark}

\section{Classification of good pairs}
In this section we will consider the odd Poisson (resp.\ generalized odd Poisson) superalgebra $PO(n,n)$ (resp.\ $PO(n,n+1)$) with the
grading of type $(0,\dots,0|1,\dots,1)$ (resp.\ $(0,\dots,0|1,\dots,1,1)$).
\begin{proposition}\label{g0spanned} Let  ${\cal P}=PO(n,n)$ or ${\cal P}=PO(n,n+1)$ and $({\cal P}, \mu)$ be a good $k$-pair.
Then  
the Lie subalgebra ${\cal P}_0$ of ${\cal P}$ is spanned by elements of the form:
$$[[\mu, a_1], \dots, a_{k-1}]b$$
with $a_1,\dots, a_{k-1}, b\in {\cal P}_{-1}$.
\end{proposition}
{\bf Proof.}  By Theorem \ref{bijection}, ${\cal P}=OP({\cal N})$ for some $k$-Nambu-Poisson algebra ${\cal N}$. Hence, by Lemma \ref{0component},  ${\cal P}_0$ is generated as a Lie algebra by elements of the form $$[[\mu, a_1], \dots, a_{k-1}]b$$
with $a_1,\dots, a_{k-1}, b\in {\cal P}_{-1}$. Let $S=\langle [[\mu, a_1], \dots, a_{k-1}] ~|~ a_1,\dots, a_{k-1}\in {\cal P}_{-1}\rangle\subset {\cal P}_0$. 

Let ${\cal P}=PO(n,n)$. Then,
for $z_1, z_2\in S$, $b_1, b_2\in {\cal P}_{-1}$, we have:
$$[z_1b_1, z_2b_2]=[z_1b_1, z_2]b_2+(-1)^{p(z_2)(p(z_1)+p(b_1)+1)}z_2[z_1b_1, b_2]=(-1)^{p(b_1)(p(z_2)+1)}[z_1, z_2]b_1b_2+$$
$$+z_1[b_1, z_2]b_2+
(-1)^{p(b_1)(p(b_2)+1)+p(z_2)(p(z_1)+p(b_1)+1)}z_2[z_1, b_2]b_1$$
since $[b_1, b_2]=0$. We recall that $[z_1, z_2]$ lies in $S$ by 
\cite[Theorem 0.2]{CantaK4}. Finally, note that $[z_1, b_2]$ and $[b_1, z_2]$ lie in ${\cal P}_{-1}$.
It follows that ${\cal P}_0\subseteq\langle[[\mu, a_1],\dots, a_{k-1}]b ~|~a_i,b\in{\cal P}_{-1}\rangle\subseteq{\cal P}_0$, hence the statement
holds for ${\cal P}=PO(n,n)$.

If ${\cal P}=PO(n,n+1)$, one uses exactly the same argument and the fact that $D_{|{\cal P}_{-1}}=0$, $D(S)\subseteq {\cal P}_{-1}$.
\hfill$\Box$

\bigskip

For any element $f\in {\cal P}_{k-1}=\F[[x_1,\dots,x_n]]\otimes\wedge^k\F^n$, we let
$f_{0}=f|_{x_1=\dots=x_n=0}\in\wedge^k\F^n$. 
 We shall say that $f$ has positive order if $f_{0}=0$. 

\begin{corollary}\label{allxi} Let ${\cal P}=PO(n,n)$ (resp.\ $PO(n,n+1)$) with the grading of type $(0,\dots,0|1,\dots,1)$
(resp.\ $(0,\dots,0|1,\dots,1,1)$).
If $\mu\in {\cal P}_{k-1}$ is such that
$\mu_0$ lies in the Grassmann subalgebra of $\wedge^k(\F^n)$ (resp.\ $\wedge^k(\F^{n+1})$) generated by some variables $\xi_{i_1}, \dots, \xi_{i_{h}}$, for
some $h<n$ (resp.\ $h<n+1$), then $\mu$ does not satisfy property G2). In particular, if $\mu_0=0$, then 
$\mu$ does not satisfy property G2).
\end{corollary}
{\bf Proof.} Suppose, on the contrary, that some $\xi_i$ does not appear in the expression of $\mu_0$. Then, by Proposition
\ref{g0spanned}, ${\cal P}_0$ does not contain $\xi_i$ and this is a contradiction since if ${\cal P}=PO(n,n)$ (resp.\
${\cal P}=PO(n,n+1)$),
 ${\cal P}_0=\langle \xi_1, \dots, \xi_n\rangle\otimes\F[[x_1, \dots, x_n]]$
(resp.\  ${\cal P}_0=\langle \xi_1, \dots, \xi_{n+1}\rangle\otimes\F[[x_1, \dots, x_{n}]]$).
\hfill$\Box$

\subsection{The case $PO(n,n)$}
In this subsection we shall determine good $k$-pairs $({\cal P}, \mu)$ for ${\cal P}=PO(n,n)$ with the $\Z_+$-grading of type $(0,\dots,0|1, \dots,1)$.
We will denote the Lie superalgebra bracket in $PO(n,n)$ simply by $[\cdot,\cdot]$.
Recall  the corresponding
description of the $\Z_+$-grading given in Example \ref{GradedGeneralized}.
When writing a monomial in $\xi_i$'s we will assume that the indices increase; elements from $\wedge^k\F^n$ will be written as linear combinations of such monomials.

\begin{lemma}\label{Step1} Let $2<k<n-1$ and
suppose  that $\mu\in PO(n,n)_{k-1}$ 
 can be written in the following form:
\begin{equation}
\mu=\xi_1\dots\xi_k+\xi_1\dots\xi_{h}\xi_{k+1}\xi_{k+2}\xi_{i_{h+1}}\dots\xi_{i_{k-2}}+\varphi+\psi,
\label{expression}
\end{equation}
where:
$$\mu_0=\xi_1\dots\xi_{k}+\xi_1\dots\xi_{h}\xi_{k+1}\xi_{k+2}\xi_{i_{h+1}}\dots\xi_{i_{k-2}}+\varphi,\,\,\varphi\in\wedge^k\F^n, \, \psi_0=0,$$ 
$h=\max\{0\leq j\leq k-2 ~|~\frac{\partial^{j+2}\mu_0}{\partial\xi_{i_1}\dots\partial\xi_{i_j}\partial\xi_{r}\partial\xi_{s}}\neq 0$,
for some $i_1<\dots<i_j\leq k,$ and some  $r,s>k\}$, 
$\frac{{\partial}^{k-1}{\varphi}}{\partial\xi_1\dots\partial\xi_{k-1}}=0$,
$\frac{{\partial}^{k}{\varphi}}{\partial\xi_1\dots\partial\xi_{h}\partial\xi_{k+1}\partial\xi_{k+2}\partial\xi_{i_{h+1}}\dots\partial\xi_{i_{k-2}}}=0$.

\bigskip

\noindent
Then $\mu$ does not satisfy property $G3)$.
\end{lemma}
{\bf Proof.} Let us first suppose that $h\geq 1$. We have:
$$[x_{k+1}, \mu]=(-1)^{h}\xi_1\dots\xi_{h}\xi_{k+2}\xi_{i_{h+1}}\dots\xi_{i_{k-2}}+\frac{\partial(\varphi+\psi)}{\partial \xi_{k+1}};$$
$$[x_{i_{k-2}}, \dots,[x_{i_{h+1}},[x_{h},\dots, [x_2, [x_1^2,[x_{k+1}, \mu]]]]]]=2(-1)^{k-2}x_1\xi_{k+2}+2x_1\frac{\partial^{k-1}(\varphi+\psi)}
{\partial \xi_{i_{k-2}}\dots\partial\xi_{i_{h+1}}\partial\xi_{h}\dots\partial\xi_1
\partial\xi_{k+1}}.$$
Therefore $[\mu,[x_{i_{k-2}}, \dots,[x_{i_{h+1}},[x_{h},\dots, [x_2, [x_1^2,[x_{k+1}, \mu]]]]]]]=$

\medskip

$=2(-1)^k(\xi_2\dots\xi_{k}((-1)^{k-2}\xi_{k+2}+
\frac{\partial^{k-1}\varphi}
{\partial \xi_{i_{k-2}}\dots\partial\xi_{i_{h+1}}\partial\xi_{h}\dots\partial\xi_1
\partial\xi_{k+1}})+$ 

\medskip

$+\xi_2\dots\xi_{h}\xi_{k+1}\xi_{k+2}\xi_{i_{h+1}}\dots\xi_{i_{k-2}}\frac{\partial^{k-1}\varphi}{\partial \xi_{i_{k-2}}\dots\partial\xi_{i_{h+1}}\partial\xi_{h}\dots\partial\xi_1
\partial\xi_{k+1}}+$

\medskip

$\frac{\partial\varphi}{\partial\xi_1}((-1)^{k-2}\xi_{k+2}+\frac{\partial^{k-1}\varphi}{\partial \xi_{i_{k-2}}\dots\partial\xi_{i_{h+1}}\partial\xi_{h}\dots\partial\xi_1
\partial\xi_{k+1}}))+\omega$,

\medskip
\noindent
for some $\omega$ of positive order.
Note that, the summand $2\xi_2\dots\xi_{k}\xi_{k+2}$
in the expression of $[\mu,[x_{i_{k-2}}, \dots,$ $[x_{i_{h+1}},[x_{h},\dots, [x_2, [x_1^2,[x_{k+1}, \mu]]]]]]]$ does not cancel out. Indeed,
due to the hypotheses on $\varphi$, the only possibility to cancel the summand $2\xi_2\dots\xi_{k}\xi_{k+2}$ is that
the expression of $\varphi$ contains the sum
$a\xi_1\dots\xi_h\xi_{k+1}\xi_{i_{h+1}}\dots\xi_{i_{k-2}}\xi_t+b\xi_1\dots\xi_{t-1}\xi_{t+1}\dots\xi_k\xi_{k+2}$, for some $t$, $2\leq t\leq k$, and some suitable coefficients $a,b\in\F^*$. But this is impossible since it is in contradiction with the maximality of $h$ if $h=k-2$, and with the hypotheses on $\varphi$ if $h<k-2$.
It follows that $[\mu,[x_{i_{k-2}}, \dots,[x_{i_{h+1}},[x_{h},\dots, [x_2, [x_1^2,[x_{k+1}, \mu]]]]]]]\neq 0$ and property $G3)$ is not satisfied.

If $h=0$, then one can use the same argument by showing that the commutator 

\noindent
$[\mu,[x_1x_{k+1},[x_{i_{1}},\dots,[x_{i_{k-2}}, \mu]]]]$ is different from zero.
\hfill$\Box$

\begin{theorem}\label{Step2}
Let ${\cal P}=PO(n,n)$. Suppose that $2<k<n-1$ and that $\mu\in PO(n,n)_{k-1}$. Then $({\cal P},\mu)$ is not a good $k$-pair. 
\end{theorem}
{\bf Proof.} 
By Corollary \ref{allxi}, if $\mu_0=0$ then $\mu$ does not satisfy property $G2)$.
Now suppose $\mu_0\neq 0$. Since $\mu_0$ lies in $\wedge^k(\F^n)$, we can assume, up to a linear change of indeterminates, that $\mu_0=\xi_1\dots\xi_{k}+f$
for some $f\in\wedge^k(\F^n)$ such that $\frac{\partial^k f}{\partial\xi_1\dots\partial\xi_k}=0$. Then, either $\mu$ does not satisfy property
$G2)$ and $({\cal P}, \mu)$ is not a good 
$k$-pair, or, again by Corollary \ref{allxi}, all $\xi_i$'s appear in the expression of $\mu_0$. Let us thus assume to be in the latter case. Then, since $k< n-1$, either
there exist some $r,s>k$ such that the indeterminates $\xi_{r}$ and $\xi_s$ both appear in the expression of $\mu_0$ in at least one monomial (case A), or
all the indeterminates $\xi_{r}$ and $\xi_s$ with $r,s>k$ appear in distinct monomials (case B).

Suppose we are in case A), and let $h=\max\{0\leq j\leq k-2 ~|~\frac{\partial^{j+2}\mu_0}{\partial\xi_{i_1}\dots\partial\xi_{i_j}\partial\xi_{r}\partial\xi_{s}}\neq 0,
 i_1<\dots<i_j\leq k; r,s>k\}$. Then we can write
$$\mu_0=\xi_1\dots\xi_k+\xi_{i_1}\dots\xi_{i_h}\xi_{r}\xi_s\xi_{i_{h+1}}\dots\xi_{i_{k-2}}+\varphi$$
for some $r,s, i_{h+1}, \dots, i_{k-2}>k$, $i_1, \dots, i_h\leq k$ and some $\varphi\in\wedge^k(\F^n)$ such that
$\frac{\partial^k \varphi}{\partial\xi_1\dots\partial\xi_k}=0$ and $\frac{\partial^k\varphi}{\partial\xi_{i_1}\dots\partial\xi_{i_h}\partial\xi_{r}\partial\xi_{s}\partial\xi_{i_1}\dots\partial\xi_{i_{k-2}}}=0$.
Up to a permutation of indices we can assume $r=k+1$, $s=k+2$, $\{i_1, \dots,i_h\}=\{1,\dots, h\}$ and up to a linear change of indeterminates we can assume
$\frac{\partial^{k-1}\varphi}{\partial\xi_{1}\dots\partial\xi_{k-1}}=0$.
Therefore $\mu$ satisfies the hypotheses  of Lemma \ref{Step1}, hence it does not satisfy property G3).

Now suppose we are in case B). Then 
$$\mu_0=\xi_1\dots\xi_{k}+\xi_{i_1}\dots\xi_{i_{k-1}}\xi_{k+1}+\xi_{j_1}\dots\xi_{j_{k-1}}\xi_{k+2}+\psi$$
for some $i_1<\dots<i_{k-1}\leq k$, $j_1<\dots<j_{k-1}\leq k$ and
$\psi\in\wedge^k(\F^n)$ such that $\frac{\partial^k \psi}{\partial\xi_1\dots\partial\xi_k}=0$,
$\frac{\partial^{k}\psi}{\partial\xi_{i_1}\dots\partial\xi_{i_{k-1}}\partial\xi_{k+1}}=0$,
$\frac{\partial^{k}\psi}{\partial\xi_{j_1}\dots\partial\xi_{j_{k-1}}\partial\xi_{k+2}}=0$,
$\frac{\partial^{2}\psi}{\partial\xi_{r}\partial\xi_{s}}=0$
for every $r,s>k$. Again by Corollary \ref{allxi}, we can assume that $\{i_1,\dots, i_{k-1}\}\neq \{j_1,\dots, j_{k-1}\}
\neq \{1,\dots, k-1\}$. Therefore there exists an index $j_l\in \{1,\dots, k\}\cap \{j_1,\dots, j_{k-1}\}$ such that $j_l\notin \{i_1,\dots, i_{k-1}\}$.

Now consider the following change of indeterminates:
$$\xi'_{j_l}=\xi_{j_l}+\xi_{k+1};~~~\xi'_j=\xi_j ~\forall j\neq j_l.$$
Then $$\mu_0=\xi'_1\dots\xi'_{k}+\xi'_{j_1}\dots\xi'_{j_{k-1}}\xi'_{k+2}+
\xi'_{j_1}\dots\xi'_{j_{i-1}}\xi'_{j_{i+1}}\dots\xi'_{j_{k-1}}\xi'_{k+1}\xi'_{k+2}+\rho$$
for some $\rho\in\wedge^k(\F^n)$ such that 
$\frac{\partial^k \rho}{\partial\xi'_1\dots\partial\xi'_k}=0$,
$\frac{\partial^{k}\rho}{\partial\xi'_{j_1}\dots\partial\xi'_{j_{k-1}}\partial\xi'_{k+2}}=0$,
$\frac{\partial^{k}\rho}{\partial\xi'_{j_1}\dots\partial\xi'_{j_{i-1}}\partial\xi'_{j_{i+1}}\dots
\partial\xi'_{k+1}\partial\xi'_{k+2}}=0$. We are now again in case A) hence the proof is concluded.
\hfill$\Box$

\begin{theorem}\label{goodHO} Let ${\cal P}=PO(n,n)$. If  $({\cal P}, \mu)$ is a good $k$-pair, then, up to isomorphisms, one of the following possibilities may occur:
\begin{itemize}
\item[a)] If $n=2h$:
\begin{itemize}
\item[a1)] $k=2$ and $\mu_0=\sum_{i=1}^h\xi_i\xi_{i+h}$;
\item[a2)] $k=n$ and $\mu_0=\xi_1\dots\xi_{n}$.
\end{itemize}
\item[b)] If $n=2h+1$:
\begin{itemize} 
\item[b1)] $k=n$ and $\mu_0=\xi_1\dots\xi_{n}$.
\end{itemize}
\end{itemize}
\end{theorem}
{\bf Proof.} By Theorem \ref{Step2}, the only possibilities for $k$ are $k=2$, $k=n-1$ or $k=n$.

By Corollary \ref{allxi}, $\frac{\partial \mu_0}{\partial\xi_i}\neq 0$ for every $i=1,\dots, n$. Using the classification of non-degenerate skew-symmetric bilinear forms, it thus follows that the case $k=2$ can occur
only if $n=2h$ and, up to equivalence, 
 $\mu_0=\sum_{i=1}^h\xi_i\xi_{i+h}$, hence we get $a1)$.

If $k=n$ then, up to rescaling the odd indeterminates, $\mu_0=\xi_1\dots\xi_n$ and we get cases $a2)$ and $b1)$.

Now assume $k=n-1$. Assume that $\frac{\partial^{n-2}\mu_0}{\partial\xi_{i_1}\dots\partial\xi_{i_{n-2}}}=\alpha\xi_{i_{n-1}}+\beta\xi_{i_n}$
for some $i_1< \dots <i_{n-2}$, $i_{n-1}<i_n$, and some $\alpha, \beta\in\F^*$. Consider the following change of indeterminates:
$$\xi'_{i_{n-1}}=\alpha\xi_{i_{n-1}}+\beta\xi_{i_n}~~~~~~~~~~\xi'_{i_j}=\xi_{i_j}\,\,\,\forall j\neq n-1.$$
Then $\frac{\partial^{n-2}\mu_0}{\partial\xi'_{i_1}\dots\partial\xi'_{i_{n-2}}}=\xi'_{i_{n-1}}$. By using induction on the lexicographic order
of the indices $i_1<\dots<i_{n-2}$, one can thus show that, up to a linear change of indeterminates, 
$\mu_0=\xi_1\dots\xi_{n-1}$, hence  $({\cal P}, \mu)$ is not a good 
$k$-pair due to
Corollary \ref{allxi}.
\hfill $\Box$

\subsection{The case $PO(n,n+1)$}
In this subsection we shall determine good pairs $({\cal P}, \mu)$ for ${\cal P}=PO(n,n+1)$ with the $\Z$-grading of type $(0,\dots,0|1, \dots,1,1)$.
We shall adopt the same notation as in the previous subsection.

\begin{lemma}\label{KStep1} Let $2\leq k<n-1$,
$\mu\in PO(n,n+1)_{k}$ and suppose that $\mu_0$ can be written in one of the following forms:
\begin{enumerate}
\item \begin{equation}
\mu_0=\xi_1\dots\xi_k\tau+\xi_1\dots\xi_{h}\xi_{k+1}\xi_{k+2}\xi_{i_{h+1}}\dots\xi_{i_{k-1}}+\varphi
\label{eq1}
\end{equation}
where:
\begin{enumerate}
\item $h=\max\{0\leq j\leq k ~|~\frac{\partial^{j+2}\mu_0}{\partial\xi_{i_1}\dots\partial\xi_{i_j}\partial\xi_{r}\partial\xi_{s}}\neq 0$, for some
 $i_1<\dots<i_j\leq k$, and $r,s>k\}$;
\item $\varphi\in\wedge^{k+1}\F^{n+1}$ is such that
$\frac{\partial^{k+1} \varphi}{\partial\xi_1\dots\partial\xi_k\partial\tau}=0$;
\end{enumerate}
\item \begin{equation}
\mu_0=\xi_1\dots\xi_k\tau+\xi_1\dots\xi_{h}\xi_{k+1}\tau\xi_{i_{h+1}}\dots\xi_{i_{k-1}}+\varphi
\label{eq2}
\end{equation}
where:
\begin{enumerate}
\item $h=\max\{0\leq j<k ~|~\frac{\partial^{j+1}\mu_0}{\partial\xi_{i_1}\dots\partial\xi_{i_j}\partial\tau}\neq 0$,
 for some $i_1<\dots<i_j\leq k\}$;
\item $\varphi\in\wedge^{k+1}\F^{n+1}$ is such that
$\frac{\partial^{k+1}\varphi}{\partial\xi_1\dots\partial\xi_h\partial\xi_{k+1}\partial\tau\partial\xi_{i_{h+1}}\dots\partial\xi_{i_{k-1}}}=0$ and
$\frac{\partial^{k+1}\varphi}{\partial\xi_2\dots\partial\xi_k\partial\tau}=0$.
\end{enumerate}
\end{enumerate}
Then $\mu$ does not satisfy property $G3)$.
\end{lemma}
{\bf Proof.} Let us first suppose that $\mu_0$ is of the form (\ref{eq1}).
Then, using the same arguments as in the proof of Lemma \ref{Step1}, one can show that
$[\mu,[x_{i_{k-1}}, \dots,[x_{i_{h+1}},[x_{h},\dots, [x_2, $ $[x_1^2,[x_{k+1}, \mu]]]]]]]\neq 0$, since in its expression the summand
$\xi_2\dots\xi_k\xi_{k+2}\tau$ does not cancel out.

Similarly, if $\mu_0$ is of the form (\ref{eq2}), then  one can show that
$[\mu,[x_{i_{k-1}}, \dots,[x_{i_{h+1}},[x_{h},\dots, [x_2, [x_1^2,$ $[1, \mu]]]]]]]\neq 0$, since in its expression the summand
$\xi_2\dots\xi_{k+1}\tau$ does not cancel out.
\hfill$\Box$

\begin{theorem}\label{KStep2}
Let ${\cal P}=PO(n,n+1)$. Suppose that $2\leq k<n-1$ and that $\mu\in {\cal P}_{k}$. Then $({\cal P},\mu)$ is not a good $(k+1)$-pair. 
\end{theorem}
{\bf Proof.} 
Let us fix a set of odd indeterminates $\xi_1, \dots, \xi_n, \xi_{n+1}=\tau$ and the corresponding basis of monomials of $\wedge(\F^{n+1})$.
By Corollary \ref{allxi}, if $\mu_0=0$ or $\frac{\partial\mu_0}{\partial\tau}=0$, then $\mu$ does not satisfy property $G2)$. 
Hence suppose that $\frac{\partial\mu_0}{\partial\tau}\neq 0$. Then we may assume, up to a linear change of indeterminates, that $\mu_0=\xi_1\dots\xi_{k}\tau+\varphi$
for some $\varphi\in\wedge^{k+1}(\F^{n+1})$ such that $\frac{\partial^{k+1} \varphi}{\partial\xi_1\dots\partial\xi_k\partial\tau}=0$. Then, either
$\frac{\partial\varphi}{\partial\tau}=0$ or $\frac{\partial\varphi}{\partial\tau}\neq 0$. 

Suppose first $\frac{\partial\varphi}{\partial\tau}=0$. Then, either for every $r,s>k$ the indeterminates $\xi_r$, $\xi_s$ appear in different monomials in the expression of $\varphi$, or there exist some $r,s>k$ such that $\xi_r$, $\xi_s$ appear in the same monomial.

In the first case $\mu_0=\xi_1\dots\xi_k\tau+\xi_1\dots\xi_k(\xi_{k+1}+\xi_{k+2})+\rho$ for some $\rho\in\wedge^{k+1}(\F^{n+1})$ such that
$\frac{\partial^{k+1} \rho}{\partial\xi_1\dots\partial\xi_k\partial\xi_{k+1}}=0=\frac{\partial^{k+1} \rho}{\partial\xi_1\dots\partial\xi_k\partial\xi_{k+2}}$.
By Corollary \ref{allxi} such an element does not satisfy property G2). Therefore we may assume that there exist some $r,s>k$ such that $\xi_r$, $\xi_s$ appear in the same monomial, i.e., that, up to a linear change of indeterminates,
$\mu_0$ is of the following form:
$$\mu_0=\xi_1\dots\xi_k\tau+\xi_1\dots\xi_h\xi_{k+1}\xi_{k+2}\xi_{i_{h+1}}\dots\xi_{i_{k-1}}+\varphi'$$
for some $\varphi'\in\wedge^{k+1}(\F^{n+1})$ such that
$\frac{\partial^{k+1}\varphi'}{\partial\xi_1\dots\partial\xi_h\partial\xi_{k+1}\partial\xi_{k+2}\partial\xi_{i_{h+1}}\dots\partial\xi_{i_{k-1}}}=0$ and
$\frac{\partial^{k+1}\varphi'}{\partial\xi_1\dots\partial\xi_k\partial\tau}=0$, where
$h=\max\{ 0\leq j\leq k ~|~\frac{\partial^{j+2}\mu_0}{\partial\xi_{i_1}\dots\partial\xi_{i_j}\partial\xi_{r}\partial\xi_{s}}\neq 0,
 \,{\mbox{for some}}\, i_1<\dots<i_j\leq k$, and $r,s>k\}$.
Therefore $\mu$ satisfies hypothesis $1.$ of Lemma \ref{KStep1}, hence it does not satisfy property G3).

Now suppose $\frac{\partial\varphi}{\partial\tau}\neq 0$. Then
$$\mu_0=\xi_1\dots\xi_k\tau+\xi_{i_1}\dots\xi_{i_h}\tau\xi_{i_{h+1}}\dots\xi_{i_k}+\psi$$
for some $i_1<\dots<i_h\leq k<i_{h+1}<\dots<i_k$, 
for some $\psi\in\wedge^{k+1}(\F^{n+1})$ such that
$\frac{\partial^{k+1}\psi}{\partial\xi_{i_1}\dots\partial\xi_{i_k}\partial\tau}=0$ and
$\frac{\partial^{k+1}\psi}{\partial\xi_1\dots\partial\xi_k\partial\tau}=0$, where
$h=\max\{0\leq j< k ~|~\frac{\partial^{j+1}\mu_0}{\partial\xi_{i_1}\dots\partial\xi_{i_j}\partial\tau}\neq 0,
 \,{\mbox{for some}}\, i_1<\dots<i_j\leq k\}$.
Now, up to a permutation of indices, we may assume that $\{i_1, \dots, i_h\}=\{1, \dots, h\}$ and $i_{h+1}=k+1$. 
Then, either $\mu$ does not satisfy property G2), or we may also assume that 
$\frac{\partial^{k}\psi}{\partial\xi_2\dots\partial\xi_k\partial\tau}=0$.
Therefore $\mu$ satisfies hypothesis $2.$ of Lemma \ref{KStep1}, hence it does not satisfy property G3).
\hfill$\Box$

\begin{theorem}\label{goodKO} Let ${\cal P}=PO(n,n+1)$. If $({\cal P}, \mu)$ is a good $(k+1)$-pair, then, up to isomorphisms, one of the following possibilities occur:
\begin{itemize}
\item[a)] If $n=2h+1$:
\begin{itemize}
\item[a1)] $k=1$ and $\mu_0=\sum_{i=1}^{h+1}\xi_i\xi_{i+h+1}$;
\item[a2)] $k=n$ and $\mu_0=\xi_1\dots\xi_{n+1}$.
\end{itemize}
\item[b)]  If $n=2h$:
\begin{itemize}
\item[b1)] $k=n$ and $\mu_0=\xi_1\dots\xi_{n+1}$.
\end{itemize}
\end{itemize}
\end{theorem}
{\bf Proof.} By Theorem \ref{KStep2}, the only possibilities for $k$ are $k=1$, $k=n-1$ or $k=n$.

By Corollary \ref{allxi}, $\frac{\partial\mu_0}{\partial\xi_i}\neq 0$ for every $i=1, \dots, n+1$. It follows that, due to the classification of non-degenerate skew-symmetric bilinear forms, the case $k=2$ can occur
only if $n=2h+1$ and, up to equivalence, 
 $\mu_0=\sum_{i=1}^{h+1}\xi_i\xi_{i+h+1}$, hence we get $a1)$.

If $k=n$ then, up to rescaling the odd indeterminates, $\mu_0=\xi_1\dots\xi_n\xi_{n+1}$ and we get cases $a2)$ and $b1)$.

Now assume $k=n-1$. Then, using the same argument as in the proof of Theorem \ref{goodHO}, one can show that, up to a linear change of indeterminates, we may assume $\mu_0=\xi_1\dots\xi_{n-1}\xi_{n+1}+f$ for some $f\in\wedge^{n}(\F^{n+1})$ such that $\frac{\partial f}{\partial\xi_{n+1}}=0$. If $f=0$ then $\mu$ does not satisfy property G2) by Corollary \ref{allxi}. If $f\neq 0$, then, up to a linear change of indeterminates, $\mu_0=\xi_1\dots\xi_{n-1}\xi_{n+1}+\xi_1\dots\xi_n=
\xi_1\dots\xi_{n-1}(\xi_{n+1}+\xi_n)$. Then, by Proposition \ref{g0spanned}, $\mu$ does not satisfy property G2). 
\hfill$\Box$

\section{The classification theorem}
\begin{remark}\label{changevariables}\em
For every invertible element $\varphi\in\F[[x_1, \dots, x_n]]$, the following  change of indeterminates preserves the odd
symplectic form, i.e., the bracket in $HO(n,n)$, and maps $\varphi\xi_1\dots\xi_n$ to $\xi'_1\dots\xi'_n$:
$$\begin{array}{lcl}
x'_1=\int_0^{x_1}\varphi^{-1}(t, x_2,\dots, x_n)dt=:\Phi, & & \xi'_1=\varphi\xi_1,\\
 &  &\\
x'_i=x_i ~~\forall~ i\neq 1, & & \xi'_i=\xi_i -\varphi\frac{\partial \Phi}{\partial x_i} \xi_1 ~~\forall~ i\neq 1.
\end{array}$$
Indeed one can check that $\{x'_i,x'_j\}_{HO}=0
=\{\xi'_i,\xi'_j\}_{HO}$ and $\{x'_i,\xi'_j\}_{HO}=\delta_{ij}$ for every $i,j=1, \dots, n$.

Note that the same change of variables, with the extra condition $\tau'=\tau$, preserves the bracket in the Lie superalgebra $KO(n,n+1)$,
and maps $\varphi\xi_1\dots\xi_n\tau$ to $\xi'_1\dots\xi'_n\tau'$. 
\end{remark}

\begin{theorem}\label{list} A complete list, up to isomorphisms, of good $k$-pairs with $k>2$, is the following:
\begin{enumerate}
\item[i)] $({\cal P}^{\varphi}, \varphi^{-1}\mu)$ with ${\cal P}=PO(n,n)$, $n>2$, $k=n$, 
$\mu=\xi_1\dots\xi_n$, $\varphi\in\F[[x_1, \dots,x_n]]$;
\item[ii)] $({\cal P}^{\varphi}, \varphi^{-1}\mu)$ with ${\cal P}=PO(n,n+1)$, $n>1$, $k=n+1$, 
$\mu=\xi_1\dots\xi_n\tau$, $\varphi\in\F[[x_1, \dots,x_n]]$. 
\end{enumerate}
\end{theorem} 
{\bf Proof.} Let ${\cal P}=PO(n,n)$ with the grading of type $(0,\dots,0|1,\dots,1)$, and let
$({\cal P}, \mu)$ be a good $k$-pair for $k>2$. 
Then, by Theorem \ref{goodHO}, we have necessarily $n>2$, $k=n$, and $\mu_0=\xi_1\dots\xi_n$. It follows that $\mu=\xi_1\dots\xi_n\psi$ for some
invertible element $\psi$ in $\F[[x_1,\dots,x_n]]$. By Remark \ref{changevariables}, up to a change of variables, we may assume $\psi=1$. 
In Example \ref{POn} we showed that
 the pair $({\cal P}, \xi_1\dots\xi_n)$ is a good $n$-pair.
Statement $i)$ then follows from Theorem \ref{rigid}, Remark \ref{gaugegradings} and Remark \ref{GaugeGood}.

Likewise,  if ${\cal P}=PO(n,n+1)$ with the grading of type $(0,\dots,0|1,\dots,1,1)$ and
$({\cal P}, \mu)$ is a good $k$-pair for $k>2$, 
by Theorem \ref{goodKO} we have necessarily $n>1$, $k=n+1$ and $\mu_0=\xi_1\dots\xi_n\tau$. It follows that $\mu=\xi_1\dots\xi_n\tau\psi$ for some
invertible element $\psi$ in $\F[[x_1,\dots,x_n]]$. Again by Remark \ref{changevariables}, we may assume $\psi=1$. Furthermore in Example \ref{POntau}
we showed that
$({\cal P}, \xi_1\dots\xi_n\tau)$ is a good $n$-pair.
Statement $ii)$ then follows from Theorem \ref{rigid}, Remark \ref{gaugegradings} and Remark \ref{GaugeGood}. \hfill$\Box$

\begin{theorem} Let $n>2$. 
\begin{itemize}
\item[a)] Any simple linearly compact generalized $n$-Nambu-Poisson algebra is gauge equivalent either to the $n$-Nambu algebra
or to the  $n$-Dzhumadildaev algebra.
\item[b)] Any simple linearly compact $n$-Nambu-Poisson algebra is isomorphic to the  $n$-Nambu algebra.
\end{itemize}
\end{theorem}
{\bf Proof.} By Theorems \ref{bijection} and \ref{rigid}, we first need to consider good $n$-pairs
$({\cal P}^{\varphi},\mu)$ where ${\cal P}=PO(k,k)$ or ${\cal P}=PO(k,k+1)$ and $n>2$. A complete list, up to isomorphisms, of such pairs
is given in Theorem \ref{list}.
The statement then follows from the construction described in Proposition \ref{reversedarrow}. We point out that 
the pair $({\cal P}^{\varphi}, \varphi^{-1}\xi_1\dots\xi_n)$, with ${\cal P}=PO(n,n)$, corresponds to ${\cal N}^{\varphi}$
where ${\cal N}$ is the $n$-Nambu algebra; similarly, the pair $({\cal P}^{\varphi}, \varphi^{-1}\xi_1\dots\xi_n\tau)$, 
with ${\cal P}=PO(n,n+1)$, corresponds to ${\cal N}^{\varphi}$,
where ${\cal N}$ is the $n$-Dzhumadildaev algebra (see also Remark \ref{correspondence}).
\hfill$\Box$

$$$$ 

\end{document}